\documentclass[12pt,a4paper,twoside,reqno]{amsart}

\title[A functorial construction of moduli of sheaves] 
{A functorial construction of\\ moduli of sheaves} 

\author{Luis \'Alvarez-C\'onsul} 
  \address{CSIC\\ Serrano 113 bis\\ 28006 Madrid, Spain} 
  \email{lac@mat.csic.es}
\author{Alastair King} 
  \address{Mathematical Sciences\\ University of Bath\\ Bath BA2 7AY, UK} 
  \email{a.d.king@maths.bath.ac.uk}

\dedicatory{Dedicated to the memory of Joseph Le Potier.}

\thanks{
The initial work was carried out at the University of
Bath, where LAC was supported by a Marie Curie Fellowship of the
European Commission.
Subsequent support has been provided by 
the European Scientific Exchange Programme of
the Royal Society of London and the Consejo Superior de
Investigaciones Cient\'{\i}ficas under Grant~15646.
LAC is partially supported by the Spanish ``Programa Ram\'on y
Cajal'' and by the Ministerio de Educaci\'on y Ciencia (Spain) under
Grant~MTM2004-07090-C03-01.
}

\usepackage{hyperref,url}
\usepackage{amssymb,latexsym}
\usepackage{amsmath,amsthm}

\usepackage[all]{xy}
\CompileMatrices
\SelectTips{cm}{12} 

\theoremstyle{plain}
\newtheorem{theorem}{Theorem}[section]
\newtheorem{lemma}[theorem]{Lemma}
\newtheorem{corollary}[theorem]{Corollary}
\newtheorem{proposition}[theorem]{Proposition}

\theoremstyle{definition}
\newtheorem{definition}[theorem]{Definition}
\newtheorem{definition-theorem}[theorem]{Definition-Theorem}

\theoremstyle{remark}
\newtheorem{remark}[theorem]{Remark}
\newtheorem*{acknowledgements}{Acknowledgements}

\newcommand{\secref}[1]{\S\ref{#1}}

\numberwithin{equation}{section}
\DeclareFontFamily{OT1}{rsfs}{}
\DeclareFontShape{OT1}{rsfs}{n}{it}{<->rsfs10}{}
\DeclareMathAlphabet{\curly}{OT1}{rsfs}{n}{it}

\newcommand{\AAA}{{\curly A}}
\newcommand{\BBB}{{\curly B}}
\newcommand{\MMM}{{\curly M}}
\newcommand{\cG}{{\mathcal G}}
\newcommand{\cK}{{\mathcal K}}
\newcommand{\cL}{{\mathcal L}}
\newcommand{\cM}{{\mathcal M}}
\newcommand{\cN}{{\mathcal N}}
\newcommand{\cO}{{\mathcal O}}
\newcommand{\cS}{{\mathcal S}}

\newcommand{\kk}{{\Bbbk}}    

\newcommand{\CC}{{\mathbb C}}
\newcommand{\PP}{{\mathbb P}}
\newcommand{\QQ}{{\mathbb Q}}

\newcommand{\FF}{{\mathbb F}}
\newcommand{\MM}{{\mathbb M}}

\newcommand{\?}{\mbox{$-$}}
\renewcommand{\(}{\left(}
\renewcommand{\)}{\right)}

\newcommand{\surj}{\to\kern-1.8ex\to}
\newcommand{\hra}{\hookrightarrow}
\newcommand{\lra}[1]{\stackrel{#1}{\longrightarrow}}
\renewcommand{\implies}{\Rightarrow} 

\newcommand{\curlyleq}{\preccurlyeq}

\newcommand{\defeq}{:=}
\newcommand{\qquot}{\operatorname{{\!/\!\!/\!}}}

\newcommand{\op}{^{\circ}} 
\newcommand{\D}{D^b}
\renewcommand{\L}{\mathbf{L}}
\newcommand{\R}{\mathbf{R}}

\newcommand{\Sch}{\operatorname{\mathbf {Sch}}}
\newcommand{\Set}{\operatorname{\mathbf {Set}}}
\newcommand{\coh}{\operatorname{\mathbf {mod}}\mbox{-}}
\newcommand{\rmod}{\operatorname{\mathbf {mod}}\mbox{-}}

\newcommand{\Hom}{\operatorname{Hom}} 
\newcommand{\End}{\operatorname{End}} 
\newcommand{\Aut}{\operatorname{Aut}} 
\newcommand{\Ext}{\operatorname{Ext}}

\newcommand{\cHom}{{\mathcal Hom}} 
\newcommand{\cEnd}{{\mathcal End}} 

\newcommand{\rk}{r} 
\newcommand{\gr}{\operatorname{gr}} 
\newcommand{\coker}{\operatorname{coker}}
\newcommand{\id}{\operatorname{id}}
\newcommand{\im}{\operatorname{im}}

\newcommand{\Spec}{\operatorname{Spec}}
\newcommand{\Proj}{\operatorname{Proj}}
\newcommand{\Supp}{\operatorname{Supp}}
\newcommand{\GL}{\operatorname{GL}}
\newcommand{\SL}{\operatorname{SL}}
\newcommand{\PGL}{\operatorname{PGL}}

\newcommand{\modfun}[2]{{\curly M}_{#2}^{#1}}
\newcommand{\modspc}[2]{{\mathcal M}_{#2}^{#1}}
\newcommand{\funpt}[1]{\underline{#1}}

\newcommand{\eval}{\varepsilon}

\newcommand{\OX}{\cO}

\newcommand{\hp}{P}

\newcommand{\brakim}[3]{{#1}^{[#3]}_{#2}}
\newcommand{\Breg}{\brakim{B}{\hp}{reg}}

\newcommand{\Rep}{R}
\newcommand{\Rreg}{\brakim{\Rep}{\hp}{reg}}

\newcommand{\Rpur}{\Rep^{pur}_\hp}

\newcommand{\Par}{Q}
\newcommand{\Pars}{\Par^{[s]}}
\newcommand{\Parss}{\Par^{[ss]}}
\newcommand{\Parpur}{\Par^{pur}}

\newcommand{\vphi}{\varphi}
\newcommand{\fmap}{f}
\newcommand{\gmap}{g}
\newcommand{\gtil}{\widetilde{\gmap}}
\newcommand{\gmapss}{g^{[ss]}}
\newcommand{\gtilss}{\widetilde{\gmap}^{[ss]}}
\newcommand{\hmap}{h}
\newcommand{\htil}{\widetilde{\hmap}}

\newcommand{\jps}{j^{\prime *}}
\newcommand{\Phiadj}{\Phi^{\vee}}

\begin{document}

\begin{abstract}
  We show how natural functors from the category of coherent sheaves on
a projective scheme to categories of Kronecker modules
can be used to construct moduli spaces of semistable sheaves.
This construction simplifies or clarifies technical 
aspects of existing constructions
and yields new simpler definitions of theta functions,
about which more complete results can be proved.
\end{abstract}

\maketitle


\section{Introduction}
\label{sec:introduction}

Let $X$ be a projective scheme over an algebraically closed
field of arbitrary characteristic.
An important set of invariants of $X$ are the projective schemes
$\modspc{ss}{X}(\hp)$, 
which are the moduli spaces for semistable coherent 
sheaves of $\cO_X$-modules with fixed Hilbert polynomial $\hp$,
with respect to a very ample invertible sheaf $\OX(1)$.

Indeed, it has been a fundamental problem to 
define and construct these moduli spaces in this generality,
ever since Mumford~\cite{Mu0} and Seshadri~\cite{Se1}
did so for smooth projective curves,
introducing the notions of stability, semistability
and S-equivalence of vector bundles.
Gieseker~\cite{Gi} and Maruyama~\cite{Ma02}
extended the definitions and constructions 
to torsion-free sheaves on higher dimensional smooth projective varieties.
Simpson~\cite{Si} completed the programme by extending 
to `pure' sheaves on arbitrary projective schemes.
Langer~\cite{La2} showed that the constructions can
be carried out in arbitrary characteristic.

Since the beginning, the method of construction has been
to identify isomorphism classes of sheaves with
orbits of a reductive group acting on a certain Quot-scheme
and then apply Geometric Invariant Theory (GIT), as developed by
Mumford~\cite{Mu1} for precisely such applications.
Thus one is required to find a projective embedding of this
Quot-scheme, with a natural linearisation of the group action,
so that the semistable sheaves correspond to GIT-semistable orbits.
One of the most natural projective embeddings to use,
as Simpson~\cite{Si} does, is into the Grassmannian
originally used by Grothendieck~\cite{Gr} to construct the Quot-scheme. 
Thus, at least in characteristic zero,
the moduli space $\modspc{ss}{X}(\hp)$ may be realised as a
closed subscheme of a GIT quotient of this Grassmannian.
(In characteristic $p$, the embedding is set-theoretic,
but not necessarily scheme-theoretic.)

The observation from which this paper grew is that the GIT quotient
of this Grassmannian has a natural moduli interpretation 
(see Remark~\ref{rem:goodquotient-Kronecker}) 
as a moduli space of (generalised) Kronecker modules or, equivalently,
modules for a certain finite dimensional associative algebra $A$.
Furthermore, the construction itself can be described in terms of a
natural functor from $\cO_X$-modules to $A$-modules. 

Taking this functorial point of view, many of the more technical aspects 
become much clearer and lead to a `one-step' construction of
the moduli spaces which is conceptually simpler than
(although structurally parallel to) 
the `two-step' process through the Quot-scheme
and the Grothendieck-Simpson embedding.

Let us begin by reviewing this two-step process.
For the first step, one chooses an integer $n$ large enough such 
that for any semistable sheaf $E$ with Hilbert polynomial $\hp$,
the natural evaluation map
\begin{equation}\label{eq:1}
  \eval_n\colon H^0(E(n)) \otimes \OX(-n) \to E
\end{equation}
is surjective and $\dim H^0(E(n))=\hp(n)$.
Thus, up to the choice of an isomorphism $H^0(E(n))\cong V$,
where $V$ is some fixed $\hp(n)$-dimensional vector space, 
we may identify $E$ with a point in the Quot-scheme parametrising
quotients of $V\otimes\OX(-n)$ with Hilbert polynomial $\hp$.
Changing the choice of isomorphism is given by the natural action
of the reductive group $\SL(V)$ on the Quot-scheme.

For the second step, one chooses another
integer $m\gg n$ so that applying the functor 
$H^0(\?\otimes\OX(m))$ to \eqref{eq:1} converts it into
a surjective map
\begin{equation}\label{eq:2}
  \alpha_E\colon H^0(E(n)) \otimes H \to H^0(E(m))
\end{equation}
where $H=H^0(\OX(m-n))$ and $\dim H^0(E(m))=\hp(m)$.
More precisely, this construction of Grothendieck's 
is applied after choosing the isomorphism $H^0(E(n))\cong V$ 
and thus $\alpha_E$ determines a point in the Grassmannian of
$\hp(m)$-dimensional quotients of $V\otimes H$.  
The fact that the Quot-scheme is embedded in the Grassmannian
means that this point determines the quotient
map $V\otimes\OX(-n) \to E$.
Now $\SL(V)$ acts
linearly on the Grassmannian and Simpson's main result \cite[Theorem
1.19]{Si} is that the GIT-semistable orbits in a certain component of
the embedded Quot-scheme are precisely those which correspond to
semistable sheaves $E$.

The essential change of strategy that 
we make in this paper is to refrain from choosing the
isomorphism $H^0(E(n))\cong V$ and 
consider just the combined effect of
the two steps above, whereby the sheaf $E$ determines functorially
the `Kronecker module' $\alpha_E$ of \eqref{eq:2}, 
from which $E$ can in turn be recovered.

A first observation in favour of this change of view-point is that 
the GIT-semistability of the orbit in the Grassmannian
is equivalent to the natural semistability of the Kronecker module $\alpha_E$
(see Remark~\ref{rem:ss-Kron-Grass}).
Thus, roughly speaking, the strategy becomes to show that
a sheaf $E$ is semistable if and only if the Kronecker module $\alpha_E$
is semistable. However, this is not quite what one is able to prove,
which indicates something of the technical difficulty 
inherent in any formulation of the construction.

Now, a Kronecker module $\alpha\colon V\otimes H\to W$
is precisely the data required to give $V\oplus W$
the structure of a (right) module for the algebra
\[
 A=\begin{pmatrix} \kk & H \\ 0 & \kk\end{pmatrix} 
\]
where $\kk$ is the base field.
If $X$ is connected and reduced, or more generally if
$H^0(\cO_X)=\kk$, then $A=\End_X(T)$, where
\[
  T=\OX(-n)\oplus\OX(-m).
\]
However, the crucial point is that $T$ is always a (left) $A$-module.
Thus, the functor
$E\mapsto \alpha_E$ is the natural functor
\begin{equation}\label{eq:HomT}
  \Phi\defeq \Phi_{n,m}=\Hom_X(T,\?)\colon \coh\cO_X\to\rmod A,
\end{equation}
where $\coh\cO_X$ is the category of coherent sheaves of
$\cO_X$-modules and $\rmod A$ is the category of finite dimensional
right $A$-modules;
see \secref{sub:Kmodules} for more discussion.

A second benefit of the functorial approach is that
the functor $\Phi$ has a left adjoint
\begin{equation}\label{eq:tensT}
  \Phiadj=\?\otimes_A T\colon \rmod A\to\coh\cO_X
\end{equation}
which provides an efficient description of how Grothendieck's
embedding works in this context.
More precisely, the fact that the $A$-module $\Hom_X(T,E)$
determines the sheaf $E$ amounts to the fact that the natural
evaluation map (the `counit' of the adjunction) 
\begin{equation}\label{eq:counit}
  \eval_E \colon  \Hom_X(T,E)\otimes_A T\to E
\end{equation}
is an isomorphism.
As we shall prove in Theorem~\ref{thm:eval=isom},
this holds not just for semistable sheaves, 
but for all $n$-regular sheaves,
in the sense of Castelnuovo-Mumford (cf. \secref{sub:MC-reg}).
Thus we see that $\Phi$ induces an `embedding' of moduli functors
\[
  \fmap_{n,m}:\modfun{reg}{X}(\hp) \to \modfun{}{A}(\hp(n),\hp(m)) 
\]
where $\modfun{reg}{X}$ is the moduli functor of
$n$-regular sheaves, of given Hilbert polynomial, 
and $\modfun{}{A}$ is the moduli functor of $A$-modules,
of given dimension vector (cf. \secref{sub:Kmodules}).
Note that semistable sheaves are $n$-regular for large enough $n$ and,
indeed, this is one condition imposed on $n$ in the usual construction.
Thus the moduli functor $\modfun{ss}{X}$ of semistable sheaves embeds in
$\modfun{reg}{X}$ and, in this way, in $\modfun{}{A}$.

Furthermore, the general machinery of adjunction provides
a simple condition for determining
when a module $M$ is in the image of the embedding:
the adjunction also has a `unit'
\begin{equation}\label{eq:unit}
 \eta_M \colon  M\to \Hom_X(T,M\otimes_A T), 
\end{equation}
and $M\cong \Hom_X(T,E)$, 
for some $E$ for which $\eval_E$ is an isomorphism,
if and only if $\eta_M$ is an isomorphism.
Hence, we can, in principle, identify which $A$-modules arise as the
image under $\Phi$ of $n$-regular sheaves and, 
in particular, can prove that the locus of such $A$-modules is 
locally closed in any family of modules
(Proposition~\ref{prop:nreg-image}).
This enables us to show that $\modfun{reg}{X}(\hp)$ is
locally isomorphic to a quotient functor (Theorem~\ref{thm:reg-loc-isomX}), 
which provides a key ingredient in our moduli space construction:
essentially replacing the Quot-scheme in the usual construction.

In particular, it is now sufficient to show that the functor $\Phi$
takes semistable sheaves to semistable $A$-modules
(one part of our main Theorem~\ref{thm:Ess-Mss}),
so that we have a locally closed embedding of functors
\[
  \fmap_{n,m}:\modfun{ss}{X}(\hp)\to \modfun{ss}{A}(\hp(n),\hp(m)) 
\]
It is known \cite{Dr,Ki} how to construct the
moduli space $\modspc{ss}{A}$ of semistable $A$-modules
in a straightforward manner
and we can then use this to construct an \emph{a priori} quasi-projective
moduli space $\modspc{ss}{X}$ of semistable sheaves
(Theorem~\ref{thm:the-theorem}).
Note that, because such moduli spaces actually parametrise S-equivalence
classes, it is helpful to observe that the functor $\Phi$
respects S-equivalence (another part of Theorem~\ref{thm:Ess-Mss}).
We complete the argument using Langton's method 
to show that $\modspc{ss}{X}$ is proper
and hence projective (Proposition~\ref{prop:MXproper}). 
Thus we obtain a closed embedding of moduli spaces
\[
  \vphi_{n,m}:\modspc{ss}{X}(\hp)\to \modspc{ss}{A}(\hp(n),\hp(m)).
\]
As a technical point, note that this embedding is scheme-theoretic 
in characteristic zero,
but in characteristic $p$ we only know that it is scheme-theoretic
on the open set of stable points (Proposition~\ref{prop:embedding-final})
and set-theoretic at the strictly semistable points.

A significant use of the embedding $\vphi$, and the whole functorial
approach, is a better understanding of theta functions, 
i.e. natural homogeneous coordinates on the moduli space.
The homogeneous coordinate rings of the moduli spaces 
$\modspc{ss}{A}$ are by now well understood 
through general results about semi-invariants of representations of
quivers \cite{DW,SV}.
More precisely, these rings are spanned by determinantal 
theta functions 
of the form
$\theta_\gamma(M)=\det\Hom_A(\gamma,M)$ for maps
\[
  \gamma\colon U_1\otimes P_1 \to U_0\otimes P_0,
\]
between projective $A$-modules.
In particular, such theta functions detect semi-stability
of $A$-modules (see Theorem~\ref{thm:SVDW}).
Now, using the adjunction 
\[
  \Hom_X(\Phiadj(\gamma),E)=\Hom_A(\gamma, \Phi(E)),
\]
we can write the restriction of such $\theta_\gamma$ to $\modspc{ss}{X}$
as an explicit theta function $\theta_\delta(E)=\det\Hom_X(\delta,E)$,
where
\[
  \delta = \Phiadj(\gamma)\colon   
  U_1\otimes\OX(-m) \to U_0\otimes\OX(-n).
\]
Thus, the theta functions $\theta_\delta$ detect semistability of sheaves
(Theorem~\ref{thm:cohochar-ssbundles}) and furthermore,
up to the same conditions as on the embedding $\vphi$,
they can be used to provide a projective embedding of $\modspc{ss}{X}$
(Theorem~\ref{thm:sep-OXmod}).

This even improves what is known about theta functions on
moduli spaces of bundles on smooth curves
(Corollary~\ref{cor:Seshadri-answer}), 
because in this case $\theta_\delta$ coincides with the usual
theta function $\theta_F$ associated to the bundle $F=\coker \delta$.

\begin{acknowledgements}
We would like to thank T.~Bridgeland, M.~Lehn,  S.~Ramanan and 
A.~Schofield for helpful remarks and expert advice.
\end{acknowledgements}

\section{Background on sheaves and Kronecker modules}
\label{sec:background}

In this section, we set out our conventions and review the notions
of semistability, stability and S-equivalence 
for sheaves (in \secref{sub:sheaves})
and for Kronecker modules (in \secref{sub:Kmodules}).  
We also explain the equivalence between $H$-Kronecker
modules and right $A$-modules.

Throughout the paper $X$ is a fixed projective scheme,
of finite type over an algebraically closed field $\kk$ of
arbitrary characteristic, with a very ample invertible sheaf $\OX(1)$. 
A `sheaf $E$ on $X$' will mean a coherent sheaf of $\cO_X$-modules.  
%

We use the notation 
\begin{align*}
  \text{``for $n\gg 0$''} &\quad \text{to mean} \quad
  \text{``$\exists n_0 \;\forall n\geq n_0$''} \quad\text{and}\\
  \text{``for $m\gg n\gg 0$''} &\quad \text{to mean} \quad
  \text{``$\exists n_0 \;\forall n\geq n_0 \;\exists m_0 \;\forall m\geq m_0$''.}
\end{align*}
Note that this notation does not necessarily imply that $n>0$ or $m>n$.

\subsection{Sheaves} 
\label{sub:sheaves}

Let $E$ be a non-zero sheaf.
Its \emph{dimension} is the dimension of the support
$
 \Supp(E):=\{ x\in X| E_x\neq 0 \}\subset X
$.
We say that $E$ is \emph{pure} if the dimension of any
non-zero subsheaf $E'\subset E$ equals the dimension of $E$.

The Hilbert polynomial $\hp(E)$ is given by
$$
 \hp(E,\ell)=\chi(E(\ell))=\sum_{i=0}^{\infty} (-1)^i h^i(E(\ell)),
$$ 
where $h^i(F)=\dim H^i(F)$.
It can be shown (e.g. \cite[Lemma 1.2.1]{HL}) that
\begin{equation}
\label{eq:31}
\hp(E,\ell)= r \ell^{d}/d! +
\textnormal{ terms of lower degree in $\ell$,}
\end{equation}
where $d$ is the dimension of $E$ and 
$r=\rk(E)$ is a positive integer,
which is roughly the `rank' of $E$,
or more strictly its
`multiplicity' \cite{HL,LP1}.

\begin{definition}\label{def:sssheaf}
A sheaf $E$ is \emph{semistable} if $E$ is
pure and, for each nonzero subsheaf $E'\subset E$,
\begin{equation}
\label{eq:54}
\frac{\hp(E')}{\rk(E')}\leq \frac{\hp(E)}{\rk(E)} 
\end{equation}
Such an $E$ is \emph{stable} if the inequality \eqref{eq:54}
is strict for all proper $E'$.
The polynomial occurring on either side of
\eqref{eq:54} is called the \emph{reduced Hilbert polynomial}
of $E'$ or $E$.
\end{definition}

In this definition, the ordering on
polynomials $p,q\in \QQ[\ell]$ is lexicographic
starting with the highest degree terms. 
Hence the inequality $p\leq q$ (resp. $p<q$) 
is equivalent to the condition 
that $p(n)\leq q(n)$ (resp. $p(n) < q(n)$) for $n\gg 0$.

Any semistable sheaf $E$ has a (not necessarily unique)
\emph{S-filtration}, 
that is, a filtration by subsheaves
\[
  0=E_0\subset E_1\subset\cdots\subset E_k=E,
\]
whose factors $E_i/E_{i-1}$ are all stable 
with the same reduced Hilbert
polynomial as $E$. 
The isomorphism class of the direct sum 
\[
  \gr E:=\bigoplus_{i=1}^k E_i/E_{i-1}
\] 
is independent of the filtration
and two semistable sheaves $E$ and $F$ are called 
\emph{S-equivalent} if $\gr E\cong \gr F$.

\begin{remark}\label{rem:Rudakov}
As observed by Rudakov \cite[\S 2]{Ru}, 
the condition of semistability may be 
formulated in a way that does not \emph{a priori} require purity.
For polynomials $p,p'$ with positive leading term, define
\begin{equation}
  \label{eq:new-ss}
p'\curlyleq p
\quad\text{iff}\quad
  \frac{p'(n)}{p'(m)} \leq \frac{p(n)}{p(m)},
  \quad\text{for $m\gg n\gg 0$.}
\end{equation}
Note: \cite{Ru} uses a different, but equivalent, 
formulation in terms of the coefficients of
the polynomials.

Then $E$ is semistable if and only if $\hp(E')\curlyleq \hp(E)$, 
for all nonzero $E'\subset E$.
In particular, a sheaf satisfying this condition is automatically pure,
because a polynomial of lower degree is bigger in this ordering.
\end{remark}

\subsection{Kronecker modules} 
\label{sub:Kmodules}

Here and throughout the paper, for integers $m>n$, 
we consider the sheaf 
\begin{equation}
  \label{eq:T}
  T=\OX(-n)\oplus\OX(-m),
\end{equation}
together with a finite dimensional $\kk$-algebra 
\begin{equation}\label{eq:A}
 A=\begin{pmatrix} \kk & H \\ 0 & \kk\end{pmatrix}
\end{equation}
of operators on $T$.
More precisely,  $A=L\oplus H\subset\End_X(T)$, where 
$L=\kk e_0\oplus\kk e_1$ is the semisimple algebra 
generated by the two projection operators onto the summands of $T$
and 
\[
  H=H^0(\OX(m-n))=\Hom(\OX(-m),\OX(-n)),
\]
acting on $T$ in the evident off-diagonal manner.
Note that $H$ is an $L$-bimodule, 
whose structure is characterised by the equation $e_0He_1=H$.
In particular, $H\otimes_L H=0$ and so 
$A$ is actually the tensor algebra of $H$
over $L$.

A right $A$-module structure on $M$ is thus
a right $L$-module structure
together with a right $L$-module map $M\otimes_L H\to M$.
The former is the same as a direct sum decomposition 
$M=V\oplus W$, where $V=Me_0$ and $W=Me_1$,
while the latter is the same as 
an $H$-Kronecker module $\alpha\colon V\otimes H\to W$.

Alternatively, we may say that $A$ is the path algebra of the quiver
\begin{equation}
\label{eq:Kron-quiv}
  \xymatrix{\bullet \ar[r]^{H} & \bullet}
\end{equation}
where $H$ is the multiplicity space for the arrow
(cf. \cite{GK}) or, after choosing a basis for $H$,
indicates that there are $\dim H$ arrows.
A representation of this quiver is precisely an $H$-Kronecker module
and the equivalence we have described is the standard one 
between representations of quivers and modules for their path algebras.

The example of particular interest is $\Hom_X(T,E)$.
On one hand, this has a natural right module structure 
over $A\subset\Hom_X(T,T)$, given by composition of maps.
On the other hand, we have the obvious decomposition
\[
 \Hom_X(T,E) = H^0(E(n))\oplus H^0(E(m))
\]
together with the multiplication map 
$\alpha_E\colon  H^0(E(n))\otimes H \to H^0(E(m))$,
as in \eqref{eq:2}.

Now, the basic discrete invariant of an $A$-module $M=V\oplus W$ 
is its \emph{dimension vector} $(\dim V, \dim W)$.
A submodule $M'\subset M$ is given by
subspaces $V'\subset V$ and $W'\subset W$ such that 
$\alpha(V'\otimes H)\subset W'$.

\begin{definition}\label{def:ssAmod}
An $A$-module $M=V\oplus W$
is \emph{semistable} if,
for each nonzero submodule $M'=V'\oplus W'$ of $M$, 
\begin{equation}
\label{eq:55}
\frac{\dim V'}{\dim W'} \leq \frac{\dim V }{\dim W }.
\end{equation}
Such a module is \emph{stable} if the inequality
is strict for all proper $M'$.
The ratio occurring on either side of \eqref{eq:55} is called 
the \emph{slope} of $M'$ or $M$. 
It lies in the (ordered) interval $[0,+\infty]$.
\end{definition}

\begin{remark}
\label{rem:ss-Kron-Grass}  
If $M$ is semistable and $\dim V>0$, 
then we see that
$\alpha\colon V\otimes H\to W$ must be surjective, because otherwise
$V\oplus\im\alpha$ is a destabilising submodule of $M$.
It is also clearly sufficient to impose the inequality~\eqref{eq:55}
just for 
\emph{saturated} submodules, 
i.e. those for which $W'=\alpha(V'\otimes H)$.
Thus, comparing Definition~\ref{def:ssAmod} 
with \cite[Proposition~1.14]{Si},
we see that (isomorphism classes of)
semistable/stable Kronecker modules correspond
precisely to GIT-semistable orbits in the Grassmannian
of $(\dim W)$-dimensional quotients of $V\otimes H$.
Of course, Definition~\ref{def:ssAmod} may be seen directly
to be equivalent to a natural GIT-semistability for Kronecker modules;
see Theorem~\ref{thm:goodquot-Kronecker} and 
Remark~\ref{rem:goodquotient-Kronecker} for further discussion.
\end{remark}

As in the case of sheaves, a semistable $A$-module
$M$ admits an S-filtration by submodules 
$$
0=M_0\subset M_1\subset\cdots\subset M_k=M,
$$
such that the quotients $M_i/M_{i-1}$ are stable 
with the same slope as $M$.
The isomorphism class of $\gr M:=\oplus_{i=1}^k 
M_i/M_{i-1}$ is independent of the filtration and 
two semistable $A$-modules $M$ and $N$ are 
\emph{S-equivalent} if $\gr M\cong \gr N$.

\begin{remark}
We shall see shortly that, for any fixed sheaf $E$,
the cohomology $H^i(E(n))$, for $i\geq1$, vanishes 
for $n\gg0$, so that the dimension vector of 
$\Hom_X(T,E)$ is then $(\hp(E,n),\hp(E,m))$.
Comparing Definition~\ref{def:ssAmod} with Remark~\ref{rem:Rudakov},
we observe that $E$ is a semistable sheaf if and only if,
for all nonzero $E'\subset E$, the $A$-submodule $\Hom_X(T,E')$
does not destabilise $\Hom_X(T,E)$, for $m\gg n\gg 0$.

This provides the basic link between semistability of sheaves
and semistability of Kronecker modules, but falls well short of what
we need. 
Our main tasks will be to show
that we can choose $n,m$ `uniformly', i.e. depending only on $\hp_E$
but not on $E$ or $E'$, and further 
to show that the submodules 
$\Hom_X(T,E')$ are the `essential' ones,
i.e. if $\Hom_X(T,E)$ is unstable, then it is destabilised by 
one of these.
\end{remark}

\begin{remark}
Kronecker modules have already played a distinguished role in the study of 
moduli of sheaves, especially on the projective plane $\PP^2$. 
Barth~\cite{Ba} showed that any stable rank 2 bundle $F$
(with $c_1=0$) could be recovered from the Kronecker module
\[
  \alpha_F\colon H^1(F(-2))\otimes H \to H^1(F(-1)),
\]
where $H=H^0(\cO(1))$,
and explicitly identified the Kronecker modules that arose
in this way. Hulek~\cite{Hu} generalised the analysis to higher rank
bundles and made an explicit link between the stability of $F$ and
the stability of $\alpha_F$.

These methods are part of the general `monad' machinery,
which was in due course used by Le Potier~\cite{LP2} to construct all
moduli spaces of sheaves on $\PP^2$, but using the more general
data of Kronecker complexes.
A different but related construction, using exceptional bundles,
enabled Drezet~\cite{Dr} to show that certain `extremal' moduli spaces
of sheaves on $\PP^2$ could be actually identified with certain 
moduli spaces of Kronecker modules.
\end{remark}

\section{The embedding functor}
\label{sec:embedding}

In this section, for $T$ and $A$ as in \eqref{eq:T} and \eqref{eq:A}
and with an additional mild assumption on $m-n$,
we show that the functor $\Hom_X(T,\?)$ embeds
sufficiently nice sheaves in the category of $A$-modules. 
Here, `sufficiently nice' means $n$-regular, in the sense of
Castelnuovo-Mumford (see \secref{sub:MC-reg}), 
and `embeds' means that $\Hom_X(T,\?)$ is fully faithful. 
In other words,
we prove, in \secref{sub:embed-reg},
that the natural evaluation map
\[
  \eval_E\colon \Hom_X(T,E)\otimes_A T\to E
\]
is an isomorphism, for any $n$-regular sheaf $E$. 
In \secref{sub:adjoint}, we describe how to construct, for any
right $A$-module $M$, the sheaf $M\otimes_A T$ and
we describe explicitly the adjunction between 
the functors $\Hom_X(T,\?)$ and $\?\otimes_A T$. 
For general background on adjoint functors, see \cite{ML}.

\subsection{Castelnuovo-Mumford regularity} 
\label{sub:MC-reg}

\begin{definition}[\cite{Mu2}]
  \label{def:regular}
  A sheaf $E$ is \emph{$n$-regular} if
  \begin{equation}
    \label{eq:regular}
     H^i(E(n-i))=0 \quad\text{for all $i>0$}. 
  \end{equation} 
 We write ``regular'' for ``$0$-regular''.
\end{definition}

Because this consists of finitely many open
conditions, it follows from Serre's Vanishing Theorem 
\cite[Theorem~III.5.2]{Ha} that any bounded family of sheaves
is $n$-regular for $n\gg 0$.
The point of this slightly odd definition
is revealed by the following consequences.

\begin{lemma}[\cite{Mu2} or {\cite[Lemma 1.7.2]{HL}}] 
\label{lem:regular}
If $E$ is $n$-regular, then
\begin{enumerate}
\item $E$ is $m$-regular for all $m\geq n$,
\item $H^i(E(n))=0$ for all $i>0$, hence $\dim H^0(E(n))=\hp(E,n)$,
\item $E(n)$ is globally generated, that is, the natural evaluation map
  $\eval_n\colon H^0(E(n))\otimes \OX(-n) \to E$
is surjective,
\item the multiplication maps
 $H^0(E(n))\otimes H^0(\OX(m-n))\to H^0(E(m))$
are surjective, for all $m\geq n$.
\end{enumerate}
\end{lemma}

In particular, by Lemma~\ref{lem:regular}(3) there is a short exact sequence
\begin{equation} \label{eq:Esyzygy}
   0\to F \lra{} V\otimes\OX(-n) \lra{\eval_n} E \to 0,
\end{equation}
where $V=H^0(E(n))$ has dimension $\hp(E,n)$.
If we also know that the `syzygy' $F$ is $m$-regular for some $m>n$,
then we can obtain a `presentation' of $E$ of the form
\begin{equation}
  \label{eq:presE}
   U\otimes\OX(-m) \lra{\delta}
   V\otimes\OX(-n) \lra{\eval_n} E \to 0,
\end{equation}
where $U=H^0(F(m))$ and $\delta$ is the composition of 
the evaluation map $f_m\colon U\otimes\OX(-m)\to F$ and the inclusion.
Note that $\hp(E)$ determines the Hilbert polynomial
of $F$ and hence the dimension of $U$.

In fact, as we shall in the proof of Theorem~\ref{thm:eval=isom}
in \secref{sub:embed-reg}, 
the map $\delta$ is essentially equivalent to the Kronecker module
$\alpha_E$ described in the Introduction
and thus we have a procedure to recover $E$ from $\alpha_E$.
This procedure can also be understood in a more functorial language
that we explain in the next subsection. 

First we note that the $m$-regularity of the syzygy $F$ 
is actually independent of $E$ and requires only that $m\gg n$.
More precisely, we have the following.

\begin{lemma} \label{lem:regker}
 Suppose $E$ is a non-zero $n$-regular sheaf,
 $F$ is the syzygy in \eqref{eq:Esyzygy} and $m>n$.
 Then $F$ is $m$-regular if and only if 
 $\OX(m-n)$ is regular.
\end{lemma}

\begin{proof}
  The fact that it is necessary for $\OX(m-n)$ to be regular
  follows immediately by applying the functor $H^i$ to the short exact sequence
\[
  0\to F(m-i)\to H^0(E(n))\otimes\OX(m-n-i) \to E(m-i)\to 0,
\]
  since $E$ is $m$-regular by Lemma~\ref{lem:regular}(1).

  On the other hand, to see that the regularity of $\OX(m-n)$ is
  sufficient, consider the following piece of the same long exact sequence
\[
  H^{i-1}(E(m-i)) \to H^i(F(m-i))\to H^0(E(n))\otimes H^{i}(\OX(m-n-i)). 
\]
  For $i>1$, the vanishing of $H^i(F(m-i))$ follows from the fact 
  $E$ is $(m-1)$-regular since $m-1\geq n$.
  In the case $i=1$, we also need that
  \[
    H^0(E(n))\otimes H^0(\OX(m-n-1))\to H^0(E(m-1))
  \]
  is surjective, which comes from Lemma~\ref{lem:regular}(4).
\end{proof}

Since we can certainly choose $m>n$ large enough that
$\OX(m-n)$ is regular, this means that every $n$-regular
sheaf $E$ with a fixed Hilbert polynomial
has a presentation by a map $\delta$ as in \eqref{eq:presE}
with fixed $U$ and $V$. 
Thus, for example, $n$-regular sheaves with given Hilbert polynomial
are bounded, since the set of such presentations certainly is.

\subsection{The adjoint functor} 
\label{sub:adjoint}

To see how to construct $M\otimes_A T$,
recall from \secref{sub:Kmodules}
that the $A$-module structure on $M$
can be specified by a direct sum decomposition $M=V\oplus W$,
giving the $L$-module structure,
together with a Kronecker module $\alpha\colon V\otimes H\to W$,
or equivalently a right $L$-module map $\alpha\colon M\otimes_L H\to M$.

On the other hand, $T$ is a left $A$-module, with its $L$-module
structure given by the decomposition $T=\OX(-n)\oplus\OX(-m)$ and
the additional $A$-module structure given by the multiplication map
\begin{equation}
  \label{eq:mu}
  \mu\colon H\otimes \OX(-m) \to \OX(-n),
\end{equation}
which we can also write as a left $L$-module map 
$\mu\colon H\otimes_L T\to T$.

From this point of view, $M\otimes_A T$ should be constructed as the 
quotient of $M\otimes_L T$ by relations
expressing the fact that the additional $H$ action is the same on
either side of the tensor product.  
More precisely, it is the cokernel of the following map.
\begin{equation}
\label{eq:10a}
  \xymatrix{
  M\otimes_L H\otimes_L T  \ar[rr]^-{1\otimes {\mu} - {\alpha}\otimes 1}
  && M\otimes_L T 
  }
\end{equation}
Writing the $L$-module structure explicitly
as a direct sum decomposition gives the following exact sequence.
\begin{equation}
\label{eq:10b}
  \xymatrix{ 
    V\otimes H\otimes \OX(-m)  \ar[rr]^-{1\otimes \mu -
      \alpha\otimes 1} && 
  }
  \begin{matrix}
    V\otimes \OX(-n)\\ \oplus \\ W\otimes \OX(-m)
  \end{matrix}
  \xymatrix{ 
    \ar[r]^-{c_n+c_m} & M\otimes_A T \to 0  
  }
\end{equation}
The exactness in \eqref{eq:10b} is also equivalent to the fact
that the following diagram is a push-out.
\begin{equation}\begin{gathered}
\label{eq:10c}
  \xymatrix{
    V\otimes \OX(-n) \ar[r]^-{c_n}
  & M\otimes_A T
 \\ V\otimes H\otimes \OX(-m) \ar[u]^{1\otimes \mu}
      \ar[r]^-{\alpha\otimes 1} 
  & W\otimes \OX(-m) \ar[u]_{c_m} 
  }
\end{gathered}\end{equation}
Thus we can see that $c_n$ is surjective if and only
if $\alpha$ is surjective.

Note that, when $\alpha$ is surjective, it carries the same information
as its kernel $\beta\colon U\to V\otimes H$, which is also
equivalent to the map 
\begin{equation}\label{eq:10d}
  \delta=(1\otimes\mu)\circ(\beta\otimes 1)
  \colon U\otimes\OX(-m) \to V\otimes\OX(-n).
\end{equation}
In this case, the kernel of $c_n$ is the image of $\delta$,
i.e. we have a presentation of $M\otimes_A T$ as the cokernel
of $\delta$, as we did for $E$ in \eqref{eq:presE}. 
Indeed, we shall see in the next section that \eqref{eq:presE}
is a special case of this construction, when $M=\Hom_X(T,E)$.

Before that, we describe explicitly the 
adjunction between $\?\otimes_A T$ and
$\Hom_X(T,\?)$,
that is, the natural isomorphism between
\[
\Hom_X(M\otimes_A T, E) \cong \Hom_A(M,\Hom_X(T,E)).
\]
This isomorphism is the first row in the
following commutative diagram with exact columns
and thus it is induced by the second and third rows.
\begin{equation}
\label{eq:adjunc}
\begin{gathered}
  \xymatrix{ 
  0 \ar[d] & 0 \ar[d]
\\
  \Hom_X(M\otimes_A T, E) \ar[r]^-{\cong} \ar[d]
& \Hom_A(M,\Hom_X(T,E)) \ar[d]
\\
  \Hom_X(M\otimes_L T, E) \ar[r]^-{\cong} \ar[d]
& \Hom_L(M,\Hom_X(T,E)) \ar[d]
\\
  \Hom_X(M\otimes_L H\otimes_L T, E) \ar[r]^-{\cong} 
& \Hom_L(M\otimes_L H,\Hom_X(T,E)) }
\end{gathered}
\end{equation}
Note that the left hand column reflects the construction of
$M\otimes_A T$ as the cokernel of \eqref{eq:10a},
while the right hand column expresses the fact that an $A$-module
map is an $L$-module map that commutes with the action of $H$.
The second and third isomorphisms and the commuting of the lower square
are clear once you unpack the $L$-module structure as a direct sum.

For example, the `unit' of the adjunction
$\eta_M\colon M\to \Hom_X(T, E))$ 
corresponds to $\id\colon E\to E$
when $E=M\otimes_A T$.
Thus
\[
  \eta_M = \eta_n\oplus\eta_m \in \Hom(V,H^0(E(n))\oplus \Hom(W,H^0(E(m))
\]
and is naturally identified with
\[
  c_n+c_m \in \Hom_X(V\otimes \OX(-n),E) \oplus \Hom_X(W\otimes
  \OX(-m),E). 
\]
On the other hand, the `counit' $\eval_E\colon M\otimes_A T \to E$,
corresponding to the identity when $M=\Hom_X(T,E)$,
is induced by the universal property of cokernels from the
map $\eval_n \oplus \eval_m\colon M\otimes_L T\to E$ given by 
the two evaluation maps
\[
  \eval_n\colon V\otimes \OX(-n)\to E
  \qquad
  \eval_m\colon W\otimes \OX(-m)\to E.
\]
Therefore, we also refer to $\eval_E$ as the evaluation map.

\subsection{Embedding regular sheaves} 
\label{sub:embed-reg}

We are now in a position to prove the main result of the section.

\begin{theorem} \label{thm:eval=isom}
 Assume that $\OX(m-n)$ is regular.
 Then the functor $\Hom_X(T,\?)$ is fully faithful on the full subcategory
 of $n$-regular sheaves.  In other words, if $E$ is an $n$-regular
 sheaf, then the natural evaluation map
 $\eval_E\colon \Hom_X(T,E)\otimes_A T\to E$ is an isomorphism.
\end{theorem}

\begin{proof}
  By Lemma~\ref{lem:regker}, the assumption means that the syzygy $F$
  in \eqref{eq:Esyzygy} is $m$-regular.
  In particular $H^1(F(m))=0$, so applying the functor
  $H^0(\?(m))$ to~\eqref{eq:Esyzygy} and recalling that $H=H^0(\OX(m-n))$,
 we obtain a short exact sequence
   \begin{equation} 
   \label{eq:4}
   \xymatrix{
     0 \ar[r] & U \ar[r]^-{\beta} & V \otimes H \ar[r]^-{\alpha} & W
     \ar[r] & 0,  
    }
    \end{equation}
    where $U=H^0(F(m))$, $V=H^0(E(n))$, $W=H^0(E(m))$
    and $\alpha$ is the Kronecker module corresponding
    to the $A$-module $\Hom_X(T,E)$.

    Now \eqref{eq:Esyzygy} and \eqref{eq:4}$\otimes\OX(-m)$
    form a commutative diagram of short
    exact sequences
    \begin{equation}
      \label{eq:cd1}
      \begin{gathered}
        \xymatrix @C=19pt {
          0\ar[r] & F \ar[r] & V\otimes \OX(-n) \ar[r]^-{\eval_n} & E \ar[r] & 0\\
          0 \ar[r] & U \otimes \OX(-m) \ar[r]^-{\beta\otimes 1}
          \ar[u]_{f_m}
          & V\otimes H\otimes \OX(-m) \ar[u]_{1\otimes \mu}
          \ar[r]^-{\alpha\otimes 1} & W\otimes \OX(-m)
          \ar[r]\ar[u]_{\eval_m} & 0 }
      \end{gathered}
    \end{equation}
    where the vertical maps are all natural evaluation maps.
    As $F$ is $m$-regular, $f_m$ is surjective and so $E$ is the cokernel
   of the map 
\[
  \delta=(1\otimes\mu)\circ(\beta\otimes 1)
\] 
of \eqref{eq:presE}, as already observed after that equation.
Alternatively, we see that the following sequence is exact.
\begin{equation}\label{eq:11}
\xymatrix{ 
 V\otimes H\otimes \OX(-m)  \ar[rr]^-{1\otimes \mu - \alpha\otimes 1} && }
 \begin{matrix}
    V\otimes \OX(-n)\\ \oplus \\ W\otimes \OX(-m)
 \end{matrix}
 \xymatrix{ \ar[r]^-{\eval_n+\eval_m} & E \ar[r] & 0  }
\end{equation}
Comparing \eqref{eq:11} with \eqref{eq:10b} we see that 
$E\cong M\otimes_A T$, where in this case $M=\Hom_X(T,E)$.
More precisely, we see that $\eval_E\colon M\otimes_A T\to E$,
as described at the end of \secref{sub:adjoint}, is an isomorphism.
\end{proof}

Note that, if $E$ is $n$-regular, then, by
Lemma~\ref{lem:regular}(1,2), the dimension vector of $\Hom_X(T,E)$ is
$(\hp(E,n),\hp(E,m))$.  Thus all $n$-regular sheaves with a fixed
Hilbert polynomial are embedded in the subcategory of $A$-modules with
a fixed dimension vector, which is a bounded subcategory:
a slight variation of the argument at the end of
\secref{sub:MC-reg}.

\begin{remark}\label{rem:gen-eval-isom}
In Theorem~\ref{thm:eval=isom}, to deduce that the evaluation map 
$\eval_E$ is an isomorphism,
we assumed that $E$ is $n$-regular and the 
syzygy $F$ in \eqref{eq:Esyzygy} is $m$-regular.
In fact, one may readily check that the proof works 
under slightly weaker hypotheses:
either that $E(n)$ is globally generated and $F$ is $m$-regular,
or that $E$ is $n$-regular and $F(m)$ is globally generated.
In the latter case, one should 
use Lemma~\ref{lem:regular}(4) to show that $\alpha$ is surjective.
\end{remark}

\section{Families and moduli}
\label{sec:families-moduli}

In this section, we show how the functorial 
embedding of the previous section
determines an embedding of moduli functors 
$f\colon \modfun{reg}{X} \to \modfun{}{A}$
(see \secref{sub:mod-fun-Kron-mod} and \secref{sub:mod-fun-sheaves}
for definitions).
This requires that we first show, in \secref{sub:presev-flat-families}, 
that Theorem~\ref{thm:eval=isom} extends to flat families
of $n$-regular sheaves.
We then show, in \secref{sub:embed-image}, how 
to identify the image of the embedding and show that it is locally closed.
We describe, in \secref{sub:mod-fun-Kron-mod}, 
how $\modfun{}{A}$ is naturally locally isomorphic
to a quotient functor of a finite dimensional vector space $R$
by a reductive group $G$ and then use \secref{sub:embed-image}
to deduce, in \secref{sub:mod-fun-sheaves}, that $\modfun{reg}{X}$
is locally isomorphic to a quotient functor of a locally closed subscheme
$Q\subset R$ by $G$. 

We conclude this section, in \secref{sub:moduli-spaces},
by describing how the GIT quotient of $R$ by $G$ is a moduli space
$\modspc{ss}{A}$ for semistable $A$-modules, corepresenting the
functor $\modfun{ss}{A}\subset \modfun{}{A}$.
This prepares the way for the construction of the moduli space of sheaves
in \secref{sec:moduli-spaces},
once we have shown that the embedding functor preserves semistability:
the main task of \secref{sec:correspondence-stability}.

\subsection{Preservation of flat families} 
\label{sub:presev-flat-families}

Let $S$ be a scheme.
A \emph{flat family $E$ over $S$ of sheaves on $X$}
is a sheaf $E$ on $X\times S$, which is flat over $S$.
On the other hand,
a \emph{flat family $M$ over $S$ of right $A$-modules}
is a sheaf $M$ of right modules over the sheaf of algebras 
$\cO_{S}\otimes A$ on $S$, which is locally free as a sheaf of
$\cO_{S}$-modules.

Let $\pi\colon X\times S\to S$ and $p_X\colon X\times S\to X$ be the
canonical projections. The adjoint pair formed by \eqref{eq:HomT}
and~\eqref{eq:tensT}, extends to an adjoint pair of functors between
the category $\rmod A\otimes\cO_S$ of sheaves of right $A$-modules on
$S$ (coherent as $\cO_{S}$-modules) 
and the category $\coh\cO_{X\times S}$ of sheaves on $X\times S$,
\begin{equation}\begin{gathered}
    \label{eq:relativeadjointpair}
    \xymatrix{ \rmod A\otimes\cO_S \ar @<-.7ex>[d]_{\?\otimes_A T}   \\
      \coh\cO_{X\times S} \ar @<-.7ex>[u]_{\Hom_X(T,\?)}  }
  \end{gathered}\end{equation}
where we are using the abbreviations
\[
  \Hom_X(T,E):=\pi_*\cHom_{X\times S}(p_X^*T,E),
\]
and
\[
  M\otimes_A T:=\pi^*M\otimes_{\cO_{X\times S}\otimes A}p_X^*T,
\] 
for a sheaf $E$ on $X\times S$ and a sheaf of right $A$-modules $M$ on
$S$.

\begin{proposition} \label{prop:fullyfaithful} 
  Assume $\OX(m-n)$ is regular. Let $S$ be any scheme.  
  Then $\Hom_X(T,\?)$ is a fully faithful functor from the
  full subcategory of $\coh\cO_{X\times S}$ consisting of flat
  families over $S$ of $n$-regular sheaves to the full subcategory of
  $\rmod A\otimes \cO_S$ consisting of flat families over $S$ of
  $A$-modules.
\end{proposition}

\begin{proof}
  To deduce that flat families are preserved by the functor, it is
  sufficient to know that $H^1(E_s(n))=0=H^1(E_s(m))$ for every sheaf
  $E_s$ in the family, so that $H^0(E_s(n))$ and $H^0(E_s(m))$ have
  locally constant dimension and hence $\Hom_X(T,E)$ is locally free.
  This vanishing follows from the regularity of $E_s(n)$, and hence
  $E_s(m)$, by Lemma~\ref{lem:regular}(1).
  
  The result then follows by applying Theorem~\ref{thm:eval=isom}
  fibrewise, using general results about cohomology and flat base
  extensions \cite{Ha}.
\end{proof}

\subsection{The image of the embedding} 
\label{sub:embed-image}

One of the uses of the functorial approach of this paper is
to identify the image of the $n$-regular sheaves by
the functor $\Hom_X(T,\?)$, using its left adjoint $\?\otimes_A T$ and
the unit of the adjunction.
Indeed, a simple general feature of adjunctions (from \cite[IV,
Theorem 1(ii)]{ML}) means that the following two statements are
equivalent:
\par\noindent
\begin{enumerate}
\item $M\cong \Hom_X(T,E)$ and
  $\eval_E \colon  \Hom_X(T,E)\otimes_A T\to E$
  is an isomorphism,
\item $E\cong M\otimes_A T$ and
  $\eta_M \colon  M\to \Hom_X(T,M\otimes_A T)$
  is an isomorphism. 
\end{enumerate}

Thus, informally speaking, any module knows how to tell that it is in
the image of the embedding and, if so, what sheaf (up to isomorphism)
it came from.  More precisely, because of
Theorem~\ref{thm:eval=isom}, a right $A$-module $M$ is isomorphic
to $\Hom_X(T,E)$ for some $n$-regular sheaf $E$ with Hilbert
polynomial $\hp$ if and only if the unit map $\eta_M$ is an
isomorphism and the sheaf $M\otimes_A T$ is $n$-regular with Hilbert
polynomial $\hp$.  It is clearly also necessary that $M$ has dimension
vector $(P(n),P(m))$.

This simple statement has an important refinement for families of
modules, which will play a key part in our subsequent construction of
the moduli of sheaves.  Roughly speaking, it says that being in the
image of the $n$-regular sheaves with Hilbert polynomial $\hp$ is a
locally closed condition in any flat family of modules.  Note that we
use the simplified notation of \secref{sub:presev-flat-families} for
applying our functors to families.

\begin{proposition}\label{prop:nreg-image}
  Assume $\OX(m-n)$ is regular.
  Let $M$ be a flat family, over a scheme $B$, of right $A$-modules
  of dimension vector $(P(n),P(m))$.
  There exists a (unique) locally closed subscheme 
  $i\colon \Breg\hra B$
  with the following properties.
  \begin{enumerate}
  \item[(a)] 
  $i^*M\otimes_A T$ is a flat family, over $\Breg$,
  of $n$-regular sheaves on $X$ with Hilbert polynomial $\hp$
  and the unit map
  \[
    \eta_{i^*M}\colon i^*M\to\Hom_X(T,i^*M\otimes_A T)
  \]
  is an isomorphism.
  \item[(b)] 
  If $\sigma\colon S\to B$ is such that $\sigma^*M\cong \Hom_X(T,E)$
  for a flat family $E$ over $S$ of $n$-regular
  sheaves on $X$ with Hilbert polynomial $\hp$,
  then $\sigma$ factors through $i\colon \Breg\hra B$
  and $E\cong \sigma^*M\otimes_A T$.
  \end{enumerate}
\end{proposition}

\begin{proof} 
  Consider $F=M\otimes_A T$, which is a sheaf on $X\times B$, but not
  necessarily flat over $B$. We can split up $B$ using the flattening
  stratification for the sheaf $F$ over the projection $X\times B\to
  B$ (see \cite{Mu2,HL} for details). Thus, there is a locally closed
  subscheme
\begin{equation}
\label{eq:flat-strat}
   j\colon B_\hp\hra B,
\end{equation}
  whose (closed) points are precisely those $b\in B$ for which the
  fibres $F_b$ have Hilbert polynomial $\hp$, and such that $j^*F$ is
  a flat family over $B_\hp$ of sheaves on $X$
  with Hilbert polynomial $\hp$.
  Therefore $B_\hp$
  contains an open set $C$ of points $b$ where the sheaf $F_b$ is
  $n$-regular and then $C$ contains an open set $D$ of points for
  which the unit map $\eta_{M_b}\colon M_b\to\Hom_X(T,F_b)$ is an
  isomorphism. The assertion that $D\subset C$ is open comes from the
  fact that, restricted to $C$, $F$ is a flat family of $n$-regular
  sheaves and so $\Hom_X(T,F)$ is a flat family of $A$-modules, by
  Proposition~\ref{prop:fullyfaithful}.  If we set $\Breg=D$, as an
  open subscheme of $B_P$, then, by construction, we have a locally
  closed subscheme $i\colon \Breg\hra B$ satisfying (a).
  
  To prove (b), first note that the counit
  \[
    \eval_E\colon \Hom_X(T,E)\otimes_A T \to E
  \]
  is an isomorphism (by Proposition~\ref{prop:fullyfaithful}) and
  hence the isomorphism $\sigma^*M\cong \Hom_X(T,E)$ implies that
  $\sigma^*M\otimes_A T \cong E$. In particular, $\sigma^*M\otimes_A T$ is a
  flat family of $n$-regular sheaves on $X$ with Hilbert polynomial
  $\hp$.  By the universal property of the flattening stratification
  (see \cite{Mu2,HL}), the fact that $\sigma^*M\otimes_A T$ is a flat
  family implies that $\sigma$ factors through $j\colon B_\hp\hra B$, while
  the fact that this is a family of $n$-regular sheaves implies that
  $\sigma$ factors through $C\subset B_\hp$.  Since the counit map for $E$
  is an isomorphism, the unit map for $\Hom_X(T,E)$ is also an
  isomorphism.  This, together with the isomorphism $\sigma^*M\cong
  \Hom_X(T,E)$, imply that the unit map 
  $\eta_{\sigma^*M}\colon \sigma^*M \to \Hom_X(T,\sigma^*M\otimes_A T)$
  is an isomorphism, hence $\sigma$
  factors through $D\subset C$.  In other words,
  $\sigma$ factors through $i\colon\Breg\hra B$, as required.
\end{proof}

\subsection{Moduli functors of Kronecker modules} 
\label{sub:mod-fun-Kron-mod}

Let $H$ be a finite dimensional vector space, 
$A$ the algebra of~\eqref{eq:A} 
and $a, b$ positive integers.  
Let $V$ and $W$ be vector spaces of dimensions $a$ and $b$, respectively.
The isomorphism classes of (right) $A$-modules, i.e. $H$-Kronecker modules,
with dimension vector $(a,b)$
are in natural bijection with the orbits of the \emph{representation space}
\begin{equation}
\label{eq:rep-space}
  \Rep \defeq \Rep_A(a,b)=\Hom_\kk(V\otimes H,W)
\end{equation}
by the canonical left action of the symmetry group
$\GL(V)\times\GL(W)$, i.e. for $g=(g_0,g_1)$ and $\alpha\in \Rep$,
\[
 g\cdot\alpha=g_1\circ \alpha \circ (g_0^{-1}\otimes 1_H).
\] 
The subgroup $\Delta=\{ (t 1, t 1) | t\in \kk^{\times}\}$
acts trivially, so we can consider the induced action of the group
\[
  G=\GL(V)\times\GL(W)/\Delta.
\]

Thus, naively, the `moduli set' parametrising
isomorphism classes of $A$-modules 
with dimension vector $(a,b)$ is the quotient set $\Rep/G$
of $G$-orbits in $\Rep$.
Unfortunately, this quotient usually cannot be given the 
geometric structure of a nice `space', e.g. a separated scheme.

After Grothendieck, an environment for giving more geometrical sense
to such quotients is the category of functors $\Sch^\circ\to\Set$
from schemes to sets.
Note first that this category does include the category of schemes itself, 
because a scheme $Z$ is determined by its functor of points
\[
  \funpt{Z}\colon \Sch^\circ\to \Set\colon X\mapsto \Hom(X,Z)
.\]
Furthermore, the Yoneda Lemma tells us that 
every natural transformation $\funpt{Y}\to\funpt{Z}$
is of the form $\funpt{f}$ for some morphism of schemes $f\colon Y\to Z$.
Note also that $\funpt{G}$ is a group valued functor,
so we may replace the quotient set $\Rep/G$ by the \emph{quotient functor}
\begin{equation}
  \label{eq:quot-fun}
  \funpt{\Rep}/\funpt{G}\colon \Sch^\circ\to \Set
  \colon X\mapsto \funpt{\Rep}(X)/\funpt{G}(X)
.\end{equation}
In this category of set-valued functors, 
we also have a replacement for the set of 
isomorphism classes of $A$-modules,
namely the \emph{moduli functor} 
\begin{equation}
\label{eq:mod-funA}
  \modfun{}{A} \defeq \modfun{}{A}(a,b) \colon \Sch\op\to \Set
,\end{equation}
where $\modfun{}{A}(S)$ is the set of
isomorphism classes of families over $S$ of $A$-modules 
with dimension vector $(a,b)$, 
in the sense of \secref{sub:presev-flat-families}. 

Note that, in this context, not all functors are as nice as the functor
of points of a scheme, which is a `sheaf' in 
an appropriate `Grothendieck topology'.
Indeed, the moduli functor $\modfun{}{A}$ 
and the quotient functor $\funpt{\Rep}/\funpt{G}$
are not strictly isomorphic, but become so after `sheafification'.

More concretely, we have the following definition
(cf. \cite[Section~1]{Si}) for the Zariski topology.

\begin{definition} \label{def:loc-isom}
  A natural transformation $g\colon \AAA\to\BBB$ between functors
  $\AAA, \BBB\colon \Sch^\circ\to\Set$ is a \emph{local isomorphism}
  if, for each $S$ in $\Sch$,
  \begin{enumerate}
  \item given $a_1,a_2\in\AAA(S)$ such that $g_S(a_1)= g_S(a_2)$,
    there is an open cover $S=\bigcup_i S_i$ such that
    $a_1|_{S_i}=a_2|_{S_i}$ for all $i$,
  \item if $b\in\BBB(S)$, then there is an open cover $S=\bigcup_i S_i$
    and $a_i\in\AAA(S_i)$ such that $g_{S_i}(a_i)=b|_{S_i}$ for all $i$.
  \end{enumerate}
\end{definition}

Since $\Rep$ carries a tautological family $\MM$ of $A$-modules
(whose fibre over each $\alpha\in \Rep$ is the $A$-module 
defined by the map $\alpha\colon V\otimes H\to W$),
we have a natural transformation
\begin{equation}\label{eq:rho}
  \hmap\colon \funpt{\Rep}\to\modfun{}{A}
\end{equation}
where $\hmap_S$ assigns to an element of $\funpt{\Rep}(S)$, 
i.e. a map $\sigma\colon S\to \Rep$, 
the isomorphism class of the pull-back $\sigma^*\MM$.

\begin{proposition}\label{prop:loc-isomA}
  The natural transformation $\hmap$ induces a local isomorphism
  $\htil\colon\funpt{\Rep}/\funpt{G} \to \modfun{}{A}$.
\end{proposition}

\begin{proof}
  First observe that two elements of $\funpt{\Rep}(S)$ define
  isomorphic families of $A$-modules if and only if they are related
  by an element of $\funpt{G}(S)$.  
  Thus the induced natural transformation is well-defined
  and satisfies part (1) of Definition~\ref{def:loc-isom}
  (even without taking open covers).  
  Part (2) is satisfied, because 
  any family of $A$-modules over $S$ can be trivialised locally,
  i.e. there is an open cover $S=\bigcup_i S_i$
  such that the restriction to each $S_i$ is
  the pull-back by a map $S_i\to \Rep$.
\end{proof}

Note that there are open subfunctors
\begin{equation}
\label{eq:sub-mod-funA}
  \modfun{s}{A} \subset \modfun{ss}{A} \subset \modfun{}{A}
\end{equation}
given by requiring that all the modules in the families over $S$ 
are stable or semistable, respectively.
Correspondingly, there are open subsets
\begin{equation}
\label{eq:sub-rep-spcA}
  \Rep^{s}\subset \Rep^{ss} \subset \Rep
\end{equation}
given by requiring that $\alpha\in \Rep$ is stable or
semistable.
Proposition~\ref{prop:loc-isomA} restricts to give local
isomorphisms
\begin{equation}
\label{eq:loc-isom-ssA}
  \funpt{\Rep^{s}}/\funpt{G} \to \modfun{s}{A}
  \qquad\text{and}\qquad
  \funpt{\Rep^{ss}}/\funpt{G} \to \modfun{ss}{A}
\end{equation}
which will be provide the route to proving,
in Theorem~\ref{thm:goodquot-Kronecker},
that there are moduli spaces of semistable and stable $A$-modules.

\subsection{Moduli functors of sheaves} 
\label{sub:mod-fun-sheaves}

Just as in \secref{sub:mod-fun-Kron-mod},
there is a straightforward moduli functor of sheaves
\begin{equation}
\label{eq:mod-fun-shfX}
  \modfun{}{X}\defeq \modfun{}{X}(P) \colon \Sch\op\to \Set
\end{equation}
which assigns to each scheme $S$ the set of isomorphism
classes of flat families over $S$ of sheaves on $X$ with
Hilbert polynomial $\hp$. 
This has open subfunctors
\begin{equation}
\label{eq:sub-mod-funX}
  \modfun{s}{X} \subset \modfun{ss}{X} \subset \modfun{}{X}
\end{equation}
given by requiring that all the sheaves in the families 
are stable or semistable, respectively.
There are also open subfunctors
\[
  \modfun{reg}{X}\defeq \modfun{reg}{X}(n) \subset\modfun{}{X}
\]
of $n$-regular sheaves, for any fixed integer $n$.

The aim of this subsection is to prove an analogue 
of Proposition~\ref{prop:loc-isomA} for $n$-regular sheaves.

\begin{theorem}\label{thm:reg-loc-isomX}
  For any $P$ and $n$, the moduli functor $\modfun{reg}{X}(P,n)$ 
  is locally isomorphic to a quotient functor $\funpt{\Par}/\funpt{G}$.
\end{theorem}

\begin{proof}
To begin, let $a=\hp(n)$, $b=\hp(m)$ and $T=\OX(-n)\oplus\OX(-m)$,
where $m>n$ and $\OX(m-n)$ is regular, 
so that Propositions~\ref{prop:fullyfaithful} and \ref{prop:nreg-image}
apply.
For $A\subset\End_X(T)$, as in~\eqref{eq:A}, consider the representation space 
$\Rep$ as in \eqref{eq:rep-space}
and let $\MM$ be the tautological family of $A$-modules on $\Rep$.

Let
\[
  \Par= \Rreg \subset \Rep
\]
be the locally closed subscheme of $\Rep$ satisfying 
Proposition~\ref{prop:nreg-image}
for the flat family $\MM$.
Roughly, $\Par$ parametrises those $\alpha\in \Rep$ that are isomorphic to
Kronecker modules arising as $\Hom_X(T,E)$ for $E$ an $n$-regular sheaf
with Hilbert polynomial $\hp$.  
More precisely, the formulation of this in 
Proposition~\ref{prop:nreg-image} gives the required result, as follows.

  Consider the following diagram of functors $\Sch^\circ\to\Set$ and
  natural transformations between them.
  \begin{equation}\label{diag:QRMM}
  \begin{gathered}
      \xymatrix{ {\funpt{\Par}} \ar[r]^{\funpt{i}}
        \ar[d]_{\gmap}&
        {\funpt{\Rep}} \ar[d]^{\hmap} \\
        \modfun{reg}{X} \ar[r]^{\fmap} & \modfun{}{A} }
  \end{gathered}
  \end{equation}
  Here $\funpt{i}$ comes from the inclusion $i\colon\Par\hra \Rep$,
  while $\hmap$ is as in \eqref{eq:rho}
  and, for any scheme $S$,
  \begin{align*}
    \gmap_S&\colon \funpt{\Par}(S) \to \modfun{reg}{X}(S)
    \colon \sigma\mapsto [\sigma^*\MM\otimes_A T], \\
    \fmap_S&\colon \modfun{reg}{X}(S)\to\modfun{}{A}(S) \colon [E] \mapsto
    [\Hom_X(T,E)],
  \end{align*}
  where $[\;]$ denotes the isomorphism class.
  
  The diagram \eqref{diag:QRMM} commutes, because,
  by Proposition~\ref{prop:nreg-image}(a),
  the unit map $\eta_{\sigma^*\MM}$ is an
  isomorphism for any $\sigma\colon S\to \Par$.
  Thus, there is a natural map
\[
  \funpt{\Par}(S) \to \modfun{reg}{X}(S)\times_{\modfun{}{A}(S)}
  \funpt{\Rep}(S) \colon \sigma\mapsto (\gmap_S(\sigma), i\circ\sigma)
\]
  Furthermore, this map is a bijection, 
  by Proposition~\ref{prop:nreg-image}(b).
  In other words, the diagram \eqref{diag:QRMM} yields 
  a pull-back in $\Set$ for all $S$, 
  i.e. it is a pull-back. 
  
  Now \eqref{diag:QRMM} induces another pull-back diagram
  \begin{equation}\label{diag:QRMM/G}
  \begin{gathered}
      \xymatrix{ {\funpt{\Par}}/{\funpt{G}} \ar[r]
        \ar[d]_{\gtil} &
        {\funpt{\Rep}}/{\funpt{G}} \ar[d]^{\htil} \\
        \modfun{reg}{X} \ar[r] & \modfun{}{A} }
  \end{gathered}
  \end{equation}
  where $\htil$ is the local
  isomorphism of Proposition~\ref{prop:loc-isomA}.
  But the pullback of a local isomorphism is a local isomorphism,
  so $\gtil$ is the local isomorphism we require.  
\end{proof}

\subsection{Moduli spaces} 
\label{sub:moduli-spaces}

It is moduli functors, as above, 
that determine what it means for a scheme to be a moduli space.
More precisely, in the terminology introduced by
Simpson~\cite[Section~1]{Si},
a moduli space is a scheme which `corepresents'
a moduli functor.

\begin{definition}\label{def:corep}
  Let $\MMM\colon \Sch^\circ\to \Set$ be a functor, $\cM$ a scheme and
  $\psi\colon \MMM\to\funpt{\cM}$ a natural transformation.
  We say that $\cM$ (or strictly $\psi$) \emph{corepresents} $\MMM$ if
  for each scheme $Y$ and each natural transformation 
  $h\colon\MMM\to\funpt{Y}$, there exists a unique
  $\sigma\colon \cM\to Y$ such that $h= \funpt{\sigma}\circ\psi$~:
\begin{equation*}
  \xymatrix{
    {\MMM} 
      \ar[rd]^-{h}
      \ar[d]_{\psi} 
    & \\
    {\funpt{\cM}} 
      \ar[r]_-{\funpt{\sigma}} &  
    {\funpt{Y}}
  }
\end{equation*}
\end{definition}

For example, suppose that an algebraic group $G$ acts on a scheme $Z$. 
Then a $G$-invariant morphism $Z\to Y$ is 
the same as a natural transformation $\funpt{Z}/ \funpt{G} \to
\funpt{Y}$, where $\funpt{Z}/ \funpt{G}$ is the quotient functor
\eqref{eq:quot-fun}.
Therefore, a $G$-invariant morphism $Z\to Z\qquot G$ is a 
`categorical quotient', 
in the sense of \cite[Definition~0.5]{Mu1},
if and only if the natural transformation $\funpt{Z}/ \funpt{G} \to
\funpt{Z\qquot G}$ corepresents the quotient functor.

The main reason for proving local isomorphism results
like Proposition~\ref{prop:loc-isomA} 
and Theorem~\ref{thm:reg-loc-isomX}
is they may be used in conjunction with the following lemma
to show that suitable categorical quotients are moduli spaces.

\begin{lemma}
\label{lem:loc-isom}
If $g\colon \AAA_1\to\AAA_2$ is a local isomorphism 
and $\psi_1\colon \AAA_1\to\funpt{Y}$ is a natural transformation, 
for a scheme $Y$, 
then there is a unique natural transformation 
$\psi_2\colon \AAA_2\to\funpt{Y}$ such that 
$\psi_1=\psi_2\circ g$~:  
\begin{equation}\begin{gathered}
\label{diag:loc-isom}
  \xymatrix{
    {\AAA_1} 
      \ar[rd]^-{\psi_1}
      \ar[d]_{g} 
    & \\
    {\AAA_2} 
      \ar[r]_-{\psi_2} &  
    {\funpt{Y}}
  }
\end{gathered}\end{equation}
Furthermore, $\psi_1$ corepresents $\AAA_1$ if and only if $\psi_2$
corepresents $\AAA_2$.
\end{lemma}

\begin{proof}
  This holds simply because `locally isomorphic'
means `isomorphic after sheafification' and $\funpt{Y}$ is a sheaf
(cf. \cite[\S 1, p. 60]{Si}).
\end{proof}

Thus, by Proposition~\ref{prop:loc-isomA}, the `moduli space' of 
all $A$-modules of dimension vector $(a,b)$ would be the categorical quotient
of $\Rep$ by $G$.
Unfortunately, this quotient collapses to just a single point, 
because every $G$ orbit in $\Rep$ has $0$ in its closure, 
so any $G$-invariant map on $\Rep$ is constant.
In fact, this collapse sensibly corresponds to the fact 
that all $A$-modules of a given dimension vector are 
Jordan-H\"older equivalent, 
i.e. they have the same simple factors.

To obtain more interesting and useful moduli spaces, 
e.g. spaces which generically 
parametrise isomorphism classes of $A$-modules,
one must restrict the class of $A$-modules that one considers.
This is typically why conditions of semistability are introduced,
so that Jordan-H\"older equivalence is replaced by S-equivalence,
which reduces to isomorphism for a larger class of objects,
namely, the stable ones rather than just the simple ones.

\begin{theorem} \label{thm:goodquot-Kronecker}
  There exist moduli spaces 
\[
  \modspc{s}{A}(a,b)\subset \modspc{ss}{A}(a,b)
\]
  of stable and semistable $A$-modules of dimension vector $(a,b)$,
  where $\modspc{ss}{A}$ is projective variety, 
  arising as a good quotient $\pi_A\colon \Rep^{ss}\to\modspc{ss}{A}$,
  and $\modspc{s}{A}$ is an open subset 
  such that the restriction $\pi_A\colon \Rep^s\to \modspc{s}{A}$
  is a geometric quotient.

  Further, the closed points of $\modspc{ss}{A}$ 
  correspond to the S-equivalence classes of semistable $A$-modules,
  and thus the closed points of $\modspc{s}{A}$ correspond
  to the isomorphism classes of stable $A$-modules.
\end{theorem}

\begin{proof}
  This is a special case of a general construction of moduli spaces
  of representations of quivers \cite{Ki}.
  A key step (\cite[Proposition~3.1]{Ki}) is
  that the open subsets $\Rep^{ss}$ and $\Rep^{s}$
  coincide with the open subsets of semistable and stable points in
  the sense of Geometric Invariant Theory~\cite{Mu1,Ne}. 
  Thus, $\modspc{ss}{A}$ can be defined as the GIT quotient of $\Rep$ by $G$.
  In particular, it is a good quotient 
  (see \cite[\S 1]{Se3} or \cite[\S 3.4]{Ne})
  of $\Rep^{ss}$ by $G$ and it is projective.
  The general machinery of GIT also ensures that 
  $\Rep^s$ is an open subset of $\Rep^{ss}$ on which the good quotient
  is geometric, yielding an open subset 
  $\modspc{s}{A}\subset \modspc{ss}{A}$.

  Furthermore, by \cite[Proposition~3.2]{Ki},
  the closed points of $\modspc{ss}{A}$ correspond to S-equivalence classes
  of semistable $A$-modules, which for stable modules are isomorphism
  classes.

  Finally, to see that $\modspc{ss}{A}$ is a moduli space in the strict
  sense,
  we note that a good quotient is, in particular, a categorical quotient
  and so we may apply Lemma~\ref{lem:loc-isom} to 
  the local isomorphism in \eqref{eq:loc-isom-ssA} to obtain
  the natural transformation $\psi_A$ in
  the following commutative diagram.
\begin{equation}\begin{gathered}
\label{diag:psiA}
  \xymatrix{
    {\funpt{\Rep^{ss}}} 
      \ar[rd]^-{\funpt{\pi_A}}
      \ar[d]_{h} 
    & \\
    {\modfun{ss}{A}} 
      \ar[r]_-{\psi_A} &  
    {\funpt{\modspc{ss}{A}}}
  }
\end{gathered}\end{equation}
The argument for $\modspc{s}{A}$ is identical.
\end{proof}

  Note that, for this special case of Kronecker modules,
  a slight variant of the construction of moduli spaces was given earlier
  by Drezet~\cite{Dr}, taking a GIT quotient of the projective space
  $\PP(\Rep)$.


\begin{remark}\label{rem:goodquotient-Kronecker}
  We claimed in the Introduction that the moduli space $\modspc{ss}{A}$ is the
  same as the GIT quotient of the Grassmannian $\cG$ of quotients of
  $V\otimes H$ of dimension $\dim(W)$ by the natural action of
  $\PGL(V)$. To justify this claim, note that there is a short exact
  sequence
  \begin{equation}
    \label{eq:ses-groups}
    1 \to \GL(W) \lra{} G \lra{} \PGL(V) \to 1.
  \end{equation}
  This enables us to make the GIT quotient of $\Rep$ by $G$ in two steps.
  In the first step, we take the GIT quotient by $\GL(W)$.  In this
  case, the semistable points $\Rep^{sur}\subset \Rep$ are also the stable
  points and are given precisely by the surjective $\alpha\colon
  V\otimes H\to W$.  Thus we have a free geometric quotient
  $\Rep^{sur}/\GL(W)$, which is the Grassmannian $\cG$ 
 (cf. \cite[\S 8.1]{Muk}).  

  The $G$-equivariant line bundle on $\Rep$ that controls
  the quotient restricts to a power of $\det W$ and thus yields a
  $\PGL(V)$-equivariant line bundle on $\cG$, which is a power of the
  $\SL(V)$-equivariant Pl\"ucker line bundle $\cO_\cG(1)$.
  Hence, the second step of taking the usual GIT quotient of $\cG$ by
  $\SL(V)$ will yield the full GIT quotient of $\Rep$ by $G$, i.e. the
  moduli space $\modspc{ss}{A}$, as required.

  Thus the GIT-semistable $\SL(V)$ orbits in $\cG$ correspond 
  to the GIT-semistable $G$ orbits in $\Rep$, which correspond to
  the semistable Kronecker modules, as Simpson effectively showed
  in~\cite[Proposition~1.14]{Si} by direct calculation
  (see also Remark~\ref{rem:ss-Kron-Grass}).
\end{remark}

\section{Preservation of semistability}
\label{sec:correspondence-stability}

In this section, we analyse the effect of the functor $\Hom_X(T,\?)$
on semistable sheaves.
We already know from \secref{sec:embedding} that the functor 
is fully faithful on $n$-regular sheaves, 
under a mild condition on $T$, i.e. on $n$ and $m$.
We start, in \secref{sub:long-list}, by listing the
much stronger conditions on $n$ and $m$ that we will now require.
All conditions have the feature that,
having fixed a Hilbert polynomial $\hp$, 
they hold for $m\gg n\gg 0$.

In \secref{sub:sheaves-modules}, we show that under these conditions
the functor $\Hom_X(T,\?)$ takes semistable sheaves of Hilbert polynomial
$\hp$ to semistable Kronecker modules and furthermore preserves 
S-filtrations.
In fact we do a little more: Theorem~\ref{thm:Ess-Mss} shows that
amongst $n$-regular pure sheaves $E$ the semistable ones are
precisely those for which $\Hom_X(T,E)$ is a semistable Kronecker module.

In \secref{sub:modules-sheaves}, we prove, for comparison,
a slightly stronger converse, replacing the assumption of $n$-regularity
by the assumption that $E=M\otimes_A T$ for some $M$.
This result is a more direct analogue of the one needed
in Simpson's construction of the moduli space and uses slightly stronger
conditions than those in \secref{sub:long-list}.

\subsection{Sufficient conditions} 
\label{sub:long-list}

For the rest of this paper, we fix a polynomial
\[
  \hp(\ell) = r\ell^d/d! + \cdots
\]
which is the Hilbert polynomial of some sheaf on $X$.
We start by noting the following variant of the 
Le Potier-Simpson estimates \cite{LP1,Si},
tailored for our later use.

\begin{theorem}
\label{thm:newLS}
There exists an integer $N_{LS}$ such that, for all $n\geq N_{LS}$
and all sheaves $E$ with Hilbert polynomial $\hp$
\begin{enumerate}
\item[(a)]
  if $E$ is pure, then the
  following conditions are equivalent:
\begin{enumerate}
\item[(1)]
  $E$ is semistable.
\item[(2)]
  $h^0(E(n))\geq \hp(n)$ and for each non-zero subsheaf $E'\subset E$,
  \[
  h^0(E'(n)) \hp \leq \hp(n) \hp(E').
  \]
\end{enumerate}
\item[(b)] 
  if $E$ is semistable and $E'$ is a non-zero subsheaf of $E$, then 
\[
  h^0(E'(n)) \hp = \hp(n) \hp(E')
    \quad\iff\quad
  \hp/r = \hp(E')/\rk(E').
\]
\end{enumerate}
\end{theorem}

\begin{proof}
  This follows from \cite[Theorem~4.4.1]{HL}, when the base field has
  characteristic zero, and \cite[Theorem~4.2]{La2} for characteristic
  $p$. 
  
  In fact, for part (a), if $E$ is semistable, then the implication
  (1)$\implies$(2) of \cite[Theorem~4.4.1]{HL} or of 
  \cite[Theorem~4.2]{La2} gives the inequality
  \begin{equation}
    \label{eq:fromLS}
      h^0(E'(n)) r \leq \hp(n) \rk(E').
  \end{equation}
  In the case when \eqref{eq:fromLS} is strict, 
  we deduce immediately that
\[
  h^0(E'(n)) \hp < \hp(n) \hp(E') 
\]
  because \eqref{eq:fromLS} gives the leading term.
  On the other hand, if we have equality in \eqref{eq:fromLS},
  then, by the last part of \cite[Theorem~4.4.1]{HL} or
  \cite[Theorem~4.2]{La2}, we have $r \hp(E') = \rk(E') \hp$, 
  and hence, with the equality in \eqref{eq:fromLS}, we get 
\[
  h^0(E'(n)) \hp = \hp(n)\hp(E').
\]

  The converse is immediate, because the leading term of the
  polynomial inequality is \eqref{eq:fromLS}, so $E$ is semistable by
  the implication (2)$\implies$(1) of \cite[Theorem~4.4.1]{HL} or of 
  \cite[Theorem~4.2]{La2}.
  
  Part (b) follows similarly from the last part of
  \cite[Theorem~4.4.1]{HL} or of \cite[Theorem~4.2]{La2}. 
\end{proof}

\begin{remark}\label{rem:dontknow}
 It is natural to ask whether the assumption that $E$ is pure
 can be dropped from the implication (2)$\implies$(1) of
 Theorem~\ref{thm:newLS}(a), or at least replaced by the assumption that $E$
 is $n$-regular.
 It is generally possible, when the dimension of $X$ is at least $3$,
 for some impure sheaf $E$ to be $n$-regular
 and yet have a subsheaf $E'$ of lower dimension with $H^0(E'(n))=0$,
 so that $E'$ does not violate condition (2).
 Thus the question is a more delicate one: having fixed $\hp(E)$,
 can $n$ be made large enough to avoid this phenomenon? 
 If so, is the condition $n\geq N_{LS}$ already sufficient?
\end{remark}

The first conditions are imposed on $n$.
\begin{enumerate}
\item[(C:1)]
  All semistable sheaves with Hilbert polynomial $\hp$ are $n$-regular.
\item[(C:2)]
  The Le Potier-Simpson estimates hold, i.e.,
  $n\geq N_{LS}$, for $N_{LS}$ as in Theorem~\ref{thm:newLS}. 
\end{enumerate}

Condition (C:1) can be satisfied because 
semistable sheaves with a fixed
Hilbert polynomial are bounded 
(see \cite{Ma1,Si,LP1} or \cite[Theorem~3.3.7]{HL} 
in characteristic zero and 
\cite[\S 4]{La1} in arbitrary characteristic),
and any bounded family of sheaves
is $n$-regular for $n\gg 0$. 
In fact, (C:1) is assumed in the proof
of the Le Potier-Simpson estimates \cite[Theorem~4.4.1]{HL}, hence of
Theorem~\ref{thm:newLS},
so (C:2) is essentially stronger than (C:1).

We next impose conditions on $m$
once $n$ has been fixed. 
The first is familiar from \secref{sec:embedding}:
\begin{enumerate}
\item[(C:3)]
 $\OX(m-n)$ is regular.
\end{enumerate}
Now, let $E$ be any $n$-regular sheaf with Hilbert polynomial $\hp$
  and 
\[
  \eval_n\colon V\otimes\OX(-n)\to E
\]
  be the evaluation map, where $V=H^0(E(n))$.
  Let $V'\subset V$ be any subspace and let
  $E'$ and $F'$ be the image and kernel of $\eval_n$ restricted to
  $V'\otimes\OX(-n)$, i.e. there is a short exact sequence
  \[
    0\to F'\to V'\otimes\OX(-n)\to E'\to 0.
  \]
  Then
\begin{enumerate}
\item[(C:4)] 
 $F'$ and $E'$ are both $m$-regular.
\item[(C:5)] 
  The polynomial relation
  \[
    \hp \dim V'\sim \hp(E')\dim V
  \]
  is equivalent to the numerical relation
  \[
    \hp(m)\dim V' \sim \hp(E',m)\dim V,
  \]
  where $\sim$ is one of $>$, $=$ or $<$.
\end{enumerate}
Condition (C:4) can be satisfied because $n$-regular sheaves $E$ 
with Hilbert polynomial $\hp$ form a bounded family 
and then so do the $V'\subset H^0(E(n))$. 
Condition (C:5) is finitely many numerical conditions on $m$ because
the set of $E'$ occurring in it is bounded and so 
there are finitely many $\hp(E')$.
Each numerical condition can be satisfied because an inequality between
polynomials in $\QQ[\ell]$ is equivalent to the same inequality with
$\ell=m$, for all sufficiently large values of $m$.

Note that (C:3) is implied by (C:4) (see Lemma~\ref{lem:regker}),
but we record it explicitly to make it clear that the results of 
\secref{sec:embedding} apply.

\subsection{From sheaves to modules} 
\label{sub:sheaves-modules}

Assuming conditions (C:1)-(C:5), we will see
how the semistability of a sheaf $E$ with Hilbert polynomial $\hp$
is related to the semistability of the $A$-module $\Hom_X(T,E)$.

Our first results explain the role of (C:4) and (C:5). 
These conditions will guarantee in particular that the `essential' subsheaves
$E'$ of an $n$-regular sheaf $E$, 
namely those with $E'(n)$ globally generated,
correspond to the `essential' submodules of $\Hom_X(T,E)$,
namely those which are `tight' in the sense of the following definition.
By ``essential'' here, we mean those which control semistability
(cf. the proofs of Propositions~\ref{prop:new-est-M} 
and \ref{prop:new-est-E} below).

\begin{definition}
  Let $M'=V'\oplus W'$ and $M''=V''\oplus W''$ be submodules
  of an $A$-module $M$. We say that $M'$ is \emph{subordinate} to $M''$ if
  \begin{equation}
    \label{eq:subord}
    V' \subset V'' \qquad\textnormal{and}\qquad W''\subset W'
  \end{equation}
  We say that $M'$ is \emph{tight} if it is subordinate to no
  submodule other than itself.
\end{definition}

Note that if $M''$ is a submodule and $V'$ and $W'$ are any subspaces
satisfying \eqref{eq:subord}, then $M'=V'\oplus W'$ is automatically a
submodule.  Furthermore, every submodule is subordinate to a tight
one and a subordinate submodule has smaller or equal slope, 
which is why the tight submodules are the `essential' ones.

\begin{lemma}\label{lem:C45}
Let $E$ be an $n$-regular sheaf with
Hilbert polynomial $\hp$ and $E'\subset E$ with $E'(n)$ globally generated.
Then 
\begin{enumerate}
\item[(a)] the evaluation map
 $\eval_{E'} \colon  \Hom_X(T,E')\otimes_A T\to E'$
 is an isomorphism, 

\item[(b)] the polynomial relation
\begin{equation}\label{eq:C45-poly<}
  h^0(E'(n)) \hp \sim \hp(n) \hp(E')
\end{equation}
  is equivalent to the numerical relation
\begin{equation}\label{eq:C45-num<}
   h^0(E'(n)) \hp(m) \sim \hp(n) h^0(E'(m)),
\end{equation}
where $\sim$ is one of $>$, $=$ or $<$,

\item[(c)] $\Hom_X(T,E')$ is a tight submodule of $\Hom_X(T,E)$. 
\end{enumerate}
\end{lemma}

\begin{proof}
  As $E$ is $n$-regular, $V=H^0(E(n))$ has dimension $\hp(n)$ and
  the evaluation map
\[
  \eval_n\colon V\otimes \OX(-n) \to E
\]
  is surjective. As $E'(n)$ is globally generated, we have a short exact 
  sequence
\begin{equation} 
  \label{eq:3'}  
  0\to F'\to V'\otimes\OX(-n)\to E'\to 0.
\end{equation}
  where $V'=H^0(E'(n))\subset V$ and the second map is the restriction
  of $\eval_n$.

For (a), note that $F'$ are $m$-regular, by (C:4). 
  Hence, by Remark~\ref{rem:gen-eval-isom},
 $\eval_{E'}$ is an isomorphism.

For (b), note that $E'$ is $m$-regular, also by (C:4).
Hence \eqref{eq:C45-num<} is equivalent to
\[
    h^0(E'(n)) \hp(m) \sim \hp(n) \hp(E',m)
\]
 which is equivalent to \eqref{eq:C45-poly<} by (C:5).

For (c), let $\alpha\colon V\otimes H\to H^0(E(m))$ be the
Kronecker module corresponding to $\Hom_X(T,E)$. 
Then, the $m$-regularity of $F'$ implies that 
$H^0(E'(m))=\alpha(V'\otimes H)$.
Thus $\Hom_X(T,E')$ is tight if there is no proper $V''\supset V'$
with $\alpha(V''\otimes H)=H^0(E'(m))$.

Suppose we have such a $V''$ and let $E''=\eval_n(V''\otimes\OX(-n))$.
As for $E'$, (C:4) implies that $E''$ is $m$-regular and
 $H^0(E''(m))= \alpha(V''\otimes H)$.
But then $E''=E'$, as both $E''(m)$ and $E'(m)$ are globally generated,
and so 
\[
  V''\subset H^0(E''(n))=H^0(E'(n))=V'.
\]
Thus $\Hom_X(T,E')$ is tight.
\end{proof}

We also have the converse of Lemma~\ref{lem:C45}(c).

\begin{lemma} \label{lem:clever}
  Let $E$ be an $n$-regular sheaf with Hilbert polynomial $\hp$ and
  $M=\Hom_X(T,E)$. If $M'\subset M$ is a tight submodule, then
  $M'=\Hom_X(T,E')$ for some subsheaf of $E'\subset E$.
  Furthermore, $E'(n)$ is globally generated and $E'\cong M'\otimes_A T$.
\end{lemma}

\begin{proof}
  Let $M'=V'\oplus W'$ be a submodule of $M=V\oplus W$,
  which corresponds to a Kronecker module $\alpha\colon V\otimes H\to W$.
  We may define a
  subsheaf $E'\subset E$ as the image of $V'\otimes\OX(-n)$ under
  the evaluation map $\eval_n\colon V\otimes\OX(-n)\to E$.
  Then $V'\subset H^0(E'(n))$ and $E'(n)$ is globally generated.

  We can apply condition (C:4) to the short exact sequence
  \[
    0\to F'\to V'\otimes\OX(-n)\to E'\to 0,
  \]
  to deduce that $F'$ is $m$-regular, and hence
  we have a short exact sequence
  \[
    0\to H^0(F'(m)) \to V'\otimes H \to H^0(E'(m)) \to 0.
  \]
  Thus $H^0(E'(m))=\alpha(V'\otimes H)\subset W'$
  and so $M'$ is subordinate to $\Hom_X(T,E')$.
  If $M'$ is tight, then $M'=\Hom_X(T,E')$, as required.  

  Even under the weaker assumption that $M'$ is saturated, 
  i.e. $W'= \alpha(V'\otimes H)$,
  we can deduce that $E'\cong M'\otimes_A T$,
  as in the proof of Theorem~\ref{thm:eval=isom}
  because $F'(m)$ is globally generated.
\end{proof}

\begin{remark}\label{rem:1-1corr}
  In Lemma~\ref{lem:clever} we could even say that $E'= M'\otimes_A T$,
  meaning that we use the natural isomorphism $M\otimes_A T \to E$
  to identify $M'\otimes_A T$
  with a subsheaf of $E$.
  Note also that, since $\?\otimes_A T$ is not generally left exact,
  it is a non-trivial fact that $M'\otimes_A T\to M\otimes_A T$
  is an injection.
  However, in this case (cf. the end of the proof) this map will actually 
  be injective whenever
  $M'$ is saturated, although then we would only have 
  $M'\subset\Hom_X(T,E')$, so such an $M'$ would not be recovered from
  $E'=M'\otimes_A T$. Indeed different saturated $M'$ could give the 
  same $E'$.

  So, as a consequence of Lemma~\ref{lem:C45}(c) and
  Lemma~\ref{lem:clever}, we can say that
  the functors $\Hom_X(T,\?)$ and $\?\otimes_A T$
  provide a one-one correspondence,
  i.e. mutually inverse bijections, 
  between the subsheaves $E'\subset E$ with $E'(n)$ globally generated
  and tight submodules $M'\subset M=\Hom_X(T,E)$.
  Because this correspondence is functorial, it automatically
  preserves inclusions between subobjects.
\end{remark}

We can now compare semistability
for $E$ and $\Hom_X(T,E)$. 

\begin{proposition}\label{prop:new-est-M}
Suppose $E$ is $n$-regular with Hilbert polynomial $\hp$.
Then $\Hom_X(T,E)$ is semistable if and only if
for all $E'\subset E$ 
\begin{equation}
\label{eq:new-est}
  h^0(E'(n)) \hp(m) \leq h^0(E'(m)) \hp(n).
\end{equation}
\end{proposition}

\begin{proof}
As $E$ is $n$-regular,
the $A$-module $M=\Hom_X(T,E)$ has dimension vector
$(\hp(n),\hp(m))$ and so this inequality is just the condition
\eqref{eq:55} for the submodule $\Hom_X(T,E')$.
Thus, if $M$ is semistable, then the inequality must hold.
On the other hand to check that $M$ is semistable, 
it is sufficient to check \eqref{eq:55} for tight 
submodules $M'$, which are all of the form $\Hom_X(T,E')$ for some $E'$,
by Lemma~\ref{lem:clever}.
\end{proof}

The next two results depend on condition (C:2),
i.e. Theorem~\ref{thm:newLS}, as well as (C:4) and (C:5).

\begin{proposition}\label{prop:new-est-E}
Suppose $E$ is $n$-regular and pure, with Hilbert polynomial $\hp$.
Then $E$ is semistable if and only if, for all $E'\subset E$,
the inequality \eqref{eq:new-est} holds.
\end{proposition}

\begin{proof}
  Suppose first that $E$ is semistable.
  To prove \eqref{eq:new-est}, we start by assuming that $E'(n)$
  is globally generated.
  By Theorem~\ref{thm:newLS}(a), we have
\begin{equation}\label{eq:fromnewLS}
  h^0(E'(n)) \hp \leq \hp(n)\hp(E'),
\end{equation}
which by Lemma~\ref{lem:C45}(b) implies \eqref{eq:new-est}.

Now for any $E'\subset E$,
let $E_0'(n)$ be the subsheaf of $E(n)$ generated by
$H^0(E'(n))$.
Then, by construction, $E_0'(n)$ is globally generated
and $H^0(E_0'(n))=H^0(E'(n))$.
Also $E_0'\subset E'$, so $h^0(E_0'(m))\leq h^0(E'(m))$,
so \eqref{eq:new-est} for $E'$ follows from 
\eqref{eq:new-est} for $E_0'$.

For the converse, we use the converse implication in
Theorem~\ref{thm:newLS}(a).
As $E$ is pure and $n$-regular, 
we have $h^0(E(n))=\hp(n)$ and 
we need to prove \eqref{eq:fromnewLS} for any $E'\subset E$.

Again, if we assume $E'(n)$ is globally generated, then
Lemma~\ref{lem:C45}(b) tells us that \eqref{eq:new-est}
implies \eqref{eq:fromnewLS}.
For arbitrary $E'\subset E$, define $E_0'$ as above.
Since we also have $\hp(E_0')\leq \hp(E')$, then \eqref{eq:fromnewLS}
for $E'$ follows from \eqref{eq:fromnewLS} for $E_0'$.
\end{proof}

\begin{proposition}\label{prop:new-ss}
Suppose $E$ is semistable with Hilbert polynomial $\hp$ and $E'\subset E$.
Then the following are equivalent:
  \begin{enumerate}
  \item[(1)] $\hp(E')/\rk(E') = \hp/r$,
  \item[(2)] $E'$ is $n$-regular and
    equality holds in \eqref{eq:new-est},
  \item[(3)] $E'(n)$ is globally generated and
    equality holds in \eqref{eq:new-est}.
  \end{enumerate}
\end{proposition}

\begin{proof}
  For (1)$\implies$(2), note that, since $E$ is semistable
  and $E'$ has the same reduced Hilbert
  polynomial, we can see that $E'$ and $E/E'$ are both semistable
  with the same reduced Hilbert polynomial.
  Hence $E'\oplus E/E'$ is semistable with Hilbert polynomial $\hp$
  and thus is $n$-regular, by (C:1).
  So in particular $E'$ is $n$-regular.
  The equality in \eqref{eq:new-est} then follows from
  the equality of reduced Hilbert polynomials.
  
  (2)$\implies$(3) immediately, 
  because regular sheaves are globally generated.

  For (3)$\implies$(1), note that $E$ is $n$-regular by (C:1).
  Since $E'(n)$ is globally generated, we can use Lemma~\ref{lem:C45}(b) 
  to deduce that equality in \eqref{eq:new-est} implies
  the polynomial equality
\[
  h^0(E'(n)) \hp = \hp(n) \hp(E'),
\]
  which implies the required equality of reduced Hilbert polynomials, by
  Theorem~\ref{thm:newLS}(b).
\end{proof}

We now combine these propositions to achieve the main aim of this section.
Note that we still assume the conditions (C:1)-(C:5).
 
\begin{theorem} \label{thm:Ess-Mss}
  Let $E$ be a sheaf on $X$ with Hilbert polynomial $\hp$.
  Then 
\begin{enumerate}
\item[(a)]
  $E$ is semistable if and only if it is $n$-regular and pure
  and the $A$-module $M=\Hom_X(T,E)$ is semistable.
\end{enumerate}
  Furthermore, if $E$ is semistable, then
\begin{enumerate}
\item[(b)] 
  $E$ is stable if and only if $M$ is stable,
\item[(c)] 
  the functors $\Hom_X(T,\?)$ and $\?\otimes_A T$
  provide a one-one correspondence
  between the subsheaves $E'\subset E$ with the same reduced Hilbert polynomial
  as $E$ and the submodules $M'\subset M$ with the same slope as $M$.
\item[(d)] 
  the functors in (c) preserve factors in this correspondence,
  i.e. if $E_1\subset E_2$ correspond to $M_1\subset M_2$, then
\begin{eqnarray*}
  M_2/M_1 &\cong & \Hom_X(T,E_2/E_1), \\
  E_2/E_1 &\cong & (M_2/M_1)\otimes_A T.
\end{eqnarray*}
\end{enumerate}
\end{theorem}

\begin{proof}
  Part (a) is the combination of
  Propositions~\ref{prop:new-est-M} and \ref{prop:new-est-E}.
  Now suppose that $E$ is semistable.
  Then (b) is just a special case of (c),
  when there are no proper subobjects on either side
  of the correspondence.

  To prove (c), first let $E'\subset E$ have the same reduced Hilbert polynomial
  as $E$. Then (1)$\implies$(3) of Proposition~\ref{prop:new-ss}
  means that $E'(n)$ is globally generated 
  and $\Hom_X(T,E')$ has the same slope as $M$.

  Conversely, if $M'\subset M$ has the same slope as $M$,
  then $M'$ is a tight submodule, because $M$ is semistable,
  and so by Lemma~\ref{lem:clever}, $M'=\Hom_X(T,E')$
  with $E'=M'\otimes_A T$ and $E'(n)$ globally generated.
  But now, by (3)$\implies$(1) of Proposition~\ref{prop:new-ss}, 
  $E'$ has the same reduced Hilbert polynomial as $E$. 

  Hence the one-one correspondence of Remark~\ref{rem:1-1corr}
  restricts to the one-one correspondence claimed here.

  To prove (d), note that the second isomorphism is automatic
  because $\?\otimes_A T$ is right exact.
  On the other hand, the stronger implication
  (1)$\implies$(2) of Proposition~\ref{prop:new-ss} tells us that
  the $E_i$ are actually $n$-regular, so
  in particular, $\Ext_X^1(T,E_1)=0$.
  Hence, when we apply $\Hom_X(T,\?)$ to the short exact sequence
  \[
   0 \to E_1 \lra{} E_2 \lra{} E_2/E_1 \to 0, 
  \]
  we obtain a short exact sequence
  \[
  0 \to  M_1 \lra{} M_2 \lra{} \Hom_X(T, E_2/E_1) \to 0, 
  \]
  which gives the first isomorphism.
\end{proof}

This theorem has the following immediate consequence.

\begin{corollary}\label{cor:grMgrE}
 Let $E$ be a semistable sheaf on $X$ with Hilbert polynomial $\hp$
 and $M=\Hom_X(T,E)$ be the corresponding semistable $A$-module. 
Then
\begin{align}
\label{eq:grM}
  \gr M &\cong \Hom_X(T,\gr E) \\
\label{eq:grE}
   \gr E &\cong (\gr M)\otimes_A T.
\end{align}
Hence a semistable sheaf $E'$ with Hilbert polynomial $\hp$
is S-equivalent to $E$ if and only if the modules $\Hom_X(T,E)$
and $\Hom_X(T,E')$ are S-equivalent.
\end{corollary}

\begin{proof}
  Let
  \begin{equation}
    \label{eq:S-filtE}
    0=E_0\subset E_1\subset\cdots\subset E_k=E
  \end{equation}
  be any S-filtration of $E$.
  Applying $\Hom_X(T,\?)$ yields a filtration
  \begin{equation}
    \label{eq:S-filtM}
    0=M_0\subset M_1\subset\cdots\subset M_k=M,
  \end{equation}
  whose terms $M_i=\Hom_X(T,E_i)$ all have the same slope as $M$.
  Indeed, this is an S-filtration of $M$, because
  if we could refine it, then applying $\?\otimes_A T$ would give
  a refinement of \eqref{eq:S-filtE}.
  Thus \eqref{eq:grM} follows, because
  $M_i/M_{i-1}=\Hom_X(T,E_i/E_{i-1})$.

  Now, \eqref{eq:grE} follows by a similar argument,
  or it follows from \eqref{eq:grM} together with
  the observation that $\gr E$ is also semistable 
  with Hilbert polynomial $\hp$ and so the counit
\[
  \Hom_X(T,\gr E) \otimes_A T \to \gr E
\]
  is an isomorphism.

  The remainder of the corollary is immediate.
\end{proof}

Thus, $\Hom_X(T,\?)$ provides a well-defined and injective set-theoretic
map from S-equivalence classes of semistable sheaves, with Hilbert
polynomial $\hp$, to S-equivalence classes of semistable modules, with
dimension vector $(\hp(n),\hp(m))$.

\begin{remark}
  It is interesting to ask whether the assumption that $E$ is
pure is really necessary (for the `if' implication) in
Proposition~\ref{prop:new-est-E} and thus in Theorem~\ref{thm:Ess-Mss}(a).
As the need for purity comes in turn from Theorem~\ref{thm:newLS}(a),
this amounts to the question asked in Remark~\ref{rem:dontknow}.

We do not know the answer, but it would perhaps be more interesting 
if purity is indeed necessary
in all these results, indicating that it might be sensible to consider
a wider notion of semistability in which certain impure sheaves could
be semistable.
\end{remark}

\subsection{From modules to sheaves} 
\label{sub:modules-sheaves}

We conclude this section with a stronger converse to
Theorem~\ref{thm:Ess-Mss}(a), 
which is included more for completeness and is not essential to
the main arguments in the paper.
The result and proof follow more closely one of the key steps in 
Simpson's construction of the moduli space 
(cf. \cite[Theorem 1.19]{Si}).

We assume conditions (C:1) and (C:3) of \secref{sub:long-list}, but
need to modify the other conditions.
First, for (C:2), we need an additional part of
\cite[Theorem~4.4.1(3)]{HL}, namely that, for all $n\geq N_{LS}$,
a sheaf $E$ with Hilbert polynomial $\hp$ is semistable, 
provided that it is pure and for each quotient
$E\surj E''$, 
\begin{equation}
\label{eq:modified-C:2}
 \frac{\hp(n)}{r}\leq\frac{h^0(E''(n))}{\rk(E'')}.
\end{equation}
We also need to strengthen conditions (C:4) and (C:5) by
requiring that they apply with $\eval_n$ replaced by any surjective
map
\[
  q\colon V\otimes\OX(-n)\to E,
\]
where $V$ is any $\hp(n)$-dimensional vector space and $E$ is any 
sheaf with Hilbert polynomial $\hp$. 
These stronger conditions can be satisfied because, once $n$ is fixed,
the sets of maps $q$ and subspaces $V'$ are bounded
(\cite{Gr}, \cite[Lemma~1.7.6]{HL}).

\begin{proposition} \label{prop:Mss-Ess}
  Let $M$ be a right $A$-module of dimension vector $(\hp(n),
  \hp(m))$ and $E$ the sheaf $M\otimes_A T$ on $X$.  If $M$ is
  semistable and $E$ is pure of Hilbert polynomial $\hp$, then $E$ is
  semistable and the unit map $\eta_M\colon M\to\Hom_X(T,E)$ is an
  isomorphism.
\end{proposition}

\begin{proof}
  Let $M=V\oplus W$, with Kronecker module $\alpha\colon V\otimes H\to W$.
  As in \secref{sub:adjoint}, let
\[
  c_n\colon V\otimes\OX(-n)\to E
  \qquad
  c_m\colon W\otimes\OX(-m)\to E
\]
  be the canonical maps, which correspond naturally to the two components
  of $\eta_M$
\[
  \eta_n\colon V\to H^0(E(n))
 \qquad
 \eta_m\colon W\to H^0(E(m)).
\]
  Since $M$ is semistable, $\alpha$ is surjective
  and so $c_n$ is surjective (see \eqref{eq:10c}).

  Let $F$ be the kernel
  \[
   0 \to F \lra{} V \otimes \OX(-n) \lra{c_n} E \to 0.  
  \]
  By the stronger version of (C:4) with $V'=V$,
   both $E$ and $F$ are $m$-regular.
  Hence, applying
  the functor $H^0(\?(m))$ to this short exact sequence,
  we see that $H^0(c_n(m))\colon V\otimes H\to H^0(E(m))$ is
  surjective and, since $H^0(c_n(m))=\eta_m\circ \alpha$, we deduce
  that $\eta_m$ is surjective.  Also $h^0(E(m))=\hp(m)=\dim W$, 
  so $\eta_m$ is an isomorphism.
  
  Hence, $\ker \eta_M=\ker \eta_n \oplus 0$ and so $\ker \eta_n=0$,
  otherwise $M$ would not be semistable.  Thus, $\eta_n$ is injective
  and to conclude that $\eta_M$ is an isomorphism, it remains to
  show that $h^0(E(n))=\hp(n)=\dim V$.
  This will follow once we have shown that $E$ is semistable, 
  and hence $n$-regular by (C:1).
  
  To prove that $E$ is semistable, we will apply the modified version
  of (C:2). 
  Let $p\colon E\surj E''$ be an epimorphism, $E'$ its kernel, and
  $V'$ and $V''$ be the kernel and the image of 
  \[
  H^0((p\circ c_n)(n))\colon V\to H^0(E''(n)), 
  \]
  respectively. 
  Then $\eta_n(V')\subset H^0(E'(n))$ and the sheaf
  \[
  E'_0=c_n(V'\otimes\OX(-n))
  \]
  is a subsheaf of $E'$. Let $c'_n\colon V'\otimes\OX(-n)\to E'_0$
  be the restriction of $c_n$ and $M'\subset M$ the submodule given by
  \[
  \alpha'=\eta_m^{-1}\circ H^0(c'_n(m))\colon V'\otimes H\to W',
  \]
  where $W'=\eta_m^{-1}(H^0(E'_0(m)))\subset W$. As $M$ is semistable,
  \[
  \frac{\dim V'}{\dim W'}\leq\frac{\dim V}{\dim W}.
  \]
  It follows from (C:4) and the fact that $\eta_m$ is an isomorphism,
  that
  \[
  \dim W'=h^0(E'_0(m))=\hp(E'_0,m),
  \]
  so the previous inequality is $\hp(m)\dim V'\leq \hp(E'_0,m) \dim
  V$. The stronger version of (C:5) 
  now implies $\hp\dim V'\leq \hp(E'_0) \dim V$.
  But $\hp(E'_0)\leq \hp(E')$, so
  \[
  \hp \dim V' \leq \hp(E')\dim V. 
  \]
  Since $\dim V = \dim V' + \dim V'' $ and $\hp=\hp(E')+\hp(E'')$,
  this is equivalent to
  \[
  \hp(E'') \dim V \leq \hp(E) \dim V''.
  \]
  If $E''$ has dimension $d$, then the leading term of this
  inequality is $\rk (E'') \dim V \leq r \dim V''$.  Since $V''
  \subset H^0(E''(n))$ and $\dim V =\hp(n)$, this implies
  \eqref{eq:modified-C:2}, so $E$ is semistable by the modified
  version of (C:2).
\end{proof}

\section{Moduli spaces of sheaves}
\label{sec:moduli-spaces}

In this section, we complete the construction of the moduli space
$\modspc{ss}{X}$ of semistable sheaves, using the formal machinery
set up in \secref{sec:families-moduli} and the results of 
\secref{sec:correspondence-stability}.
In \secref{sub:construction-moduli}, we use the moduli space
$\modspc{ss}{A}$ of semistable $A$-modules, 
as described in \secref{sub:moduli-spaces},
to construct the moduli space $\modspc{ss}{X}$.
In characteristic zero, this is simply a locally closed subscheme
of $\modspc{ss}{A}$, but in characteristic $p$ the construction
is a little more delicate. 
At this stage we only know that $\modspc{ss}{X}$ is a quasi-projective 
scheme. 
In \secref{sub:proper}, we use Langton's method
to show that $\modspc{ss}{X}$ is proper and hence projective.
In \secref{sub:conc-emb}, we look more closely at what we can say about the
embedding of moduli spaces $\vphi\colon \modspc{ss}{X} \to \modspc{ss}{A}$,
in particular when it is a scheme-theoretic embedding and when
only set-theoretic.
In \secref{sub:uni-properties} and \secref{sub:relative-moduli}, 
we discuss some technical enhancements to the construction.

\subsection{Construction of the moduli spaces} 
\label{sub:construction-moduli}

 
To construct the moduli space of
semistable sheaves on $X$, 
we start from the analogue of Theorem~\ref{thm:reg-loc-isomX}
for semistable and stable sheaves.

\begin{proposition}\label{prop:sss-loc-isomX}
  For any $P$, the moduli functors
 $\modfun{s}{X}(P)\subset \modfun{ss}{X}(P)$
  are locally isomorphic to quotient functors
 $\funpt{\Pars}/\funpt{G}\subset \funpt{\Parss}/\funpt{G}$.
\end{proposition}

\begin{proof}
  Choose $n$ large enough that all semistable sheaves of Hilbert
  polynomial $\hp$ are $n$-regular and choose $m$ large enough
  that $\OX(m-n)$ is regular.
  Define open subsets
\[
  \Pars\subset \Parss\subset \Par,
\]
with $\Par$ as in Theorem~\ref{thm:reg-loc-isomX},
to be the open loci where the fibres
of the tautological flat family
$\FF=i^*\MM\otimes_A T$ over $\Par$,
are stable and semistable, respectively.

Then the restriction of the family $\FF$ to $\Parss$ 
determines a natural transformation 
\begin{equation}
\label{eq:gmap-ss}
  \gmapss\colon \funpt{\Parss} \to \modfun{ss}{X}
,\end{equation}
which is the restriction of the natural transformation 
$\gmap\colon \funpt{\Par} \to \modfun{reg}{X}$ of \eqref{diag:QRMM}
and, by definition of $\Parss$, is also the pull-back of $\gmap$
along $\modfun{ss}{X}\hra \modfun{reg}{X}$.
Hence we obtain a local isomorphism
\begin{equation}
  \label{eq:gtil-ss}
  \gtilss\colon\funpt{\Parss}/\funpt{G} \to \modfun{ss}{X},
\end{equation}
by restricting the local isomorphism 
$\gtil\colon \funpt{\Par}/\funpt{G} \to \modfun{reg}{X}$ 
of \eqref{diag:QRMM/G}.

An identical argument applies to $\Pars$ and $\modfun{s}{X}$.
\end{proof}

Thus, the main step in constructing the moduli space is 
to show that $\Parss$ has a good quotient
and, further, that this 
restricts to a geometric quotient of $\Pars$.
Before we do this, we note the following (well-known) lemma.

\begin{lemma}\label{lem:goodquotient}
  Let $\pi\colon Z\to Z\qquot G$ be a good quotient for the action
  of a reductive algebraic group $G$ on a scheme $Z$
  and let $Y$ be a $G$-invariant open subset of $Z$.
  Suppose further that,
  for each $G$-orbit $O$ in Y, the closed orbit $O'\subset Z$
  contained in the orbit closure $\overline{O}$ is also in Y.
  Then $\pi$ restricts to a good quotient $Y\to Y\qquot G$,
  where $Y\qquot G=\pi(Y)$ is an open subset of $Z\qquot G$.
\end{lemma}
               
\begin{proof}
  For any good quotient, $\pi(Z\setminus Y)$ will be closed in $Z\qquot G$. 
  Since $\pi$ induces a bijection between the closed
  orbits in $Z$ and the (closed) points in $Z\qquot G$ and $\pi(O)=\pi(O')$,
  the additional assumption here implies that $\pi(Y)$ and
  $\pi(Z\setminus Y)$ are disjoint. 
  As $\pi$ is surjective,
\[
  \pi(Y) = (Z\qquot G)\setminus\pi(Z\setminus Y)
\]
  and so $\pi(Y)$ is open in $Z\qquot G$.
  Furthermore, $Y=\pi^{-1}(\pi(Y))$ and so 
  the restriction $\pi\colon Y\to \pi(Y)$ is a good quotient,
  because this is a property which is local in $Z\qquot G$
  (cf. \cite[Definition 1.5]{Se3}).
\end{proof}

To apply this lemma, we now suppose that $n,m$ satisfy the conditions
(C:1)-(C:5) of \secref{sub:long-list}, so $\Parss$ is locally closed
in $\Rep^{ss}$, by Theorem~\ref{thm:Ess-Mss}.

\begin{proposition}\label{prop:goodquotient}
  The good quotient $\pi_A\colon \Rep^{ss}\to\modspc{ss}{A}$ of
  Theorem~\ref{thm:goodquot-Kronecker}
  and the inclusion $i\colon \Parss\to \Rep^{ss}$ determine
  a (unique) commutative diagram
\begin{equation}\begin{gathered}
\label{diag:goodquot}
  \xymatrix{ 
    \Parss \ar[r]^{i} \ar[d]_{\pi_X}
 &  {\Rep}^{ss} \ar[d]^{\pi_A}
 \\ \modspc{ss}{X} \ar[r]^{\vphi} 
 &  \modspc{ss}{A} 
  }
\end{gathered}\end{equation}
  where $\modspc{ss}{X}$ is quasi-projective,
  $\pi_X$ is a good quotient and
  $\vphi$ induces a set-theoretic injection on closed points.
  In characteristic zero, $\vphi$ is the
  inclusion of the locally closed subscheme $\pi_A\(\Parss\)$.
\end{proposition}

\begin{proof}
  Let $Y=\Parss$ and $\overline{Y}$ be the closure of $Y$ in $\Rep$.
  Let $Z=\overline{Y}\cap \Rep^{ss}$, which is closed in $\Rep^{ss}$.
 Observe that the additional assumption of Lemma~\ref{lem:goodquotient}
 holds in this case.
 The closed orbit in the orbit closure of a
 point in $\Parss$, corresponding to a module $M=\Hom_X(T,E)$, is the
 orbit corresponding to the associated graded module $\gr M$ 
 (by \cite[Proposition~3.2]{Ki}).  
 But $\gr M \cong \Hom_X(T,\gr E)$,
 by Corollary~\ref{cor:grMgrE},
 and $\gr E$ is semistable,
 so this closed orbit is also in $\Parss$.

 Thus, to obtain the quasi-projective good quotient
 $\modspc{ss}{X}=Y\qquot G$ with the given additional properties,
 it is sufficient, using Lemma~\ref{lem:goodquotient}, to 
 show that the closed subscheme $Z\subset \Rep^{ss}$ has a 
 projective good quotient with corresponding additional properties.

 In characteristic zero, the processes of taking quotient rings 
 and invariant subrings commute (as seen by using the Reynolds operator),
 and so the (scheme-theoretic) image $\pi_A(Z)$ is the good quotient
 $Z\qquot G$.
 Furthermore, $\pi_A(Z)$ is a closed subscheme of $\modspc{ss}{A}$,
 which is projective, and hence $Z\qquot G$ is projective.

 Thus we have a commutative diagram
  \begin{equation}\label{diag:goodquotbar}
  \begin{gathered}
      \xymatrix{ Z \ar[r]^{\overline{i}} \ar[d]_{\pi} &
        \Rep^{ss} \ar[d]^{\pi_A} \\
        Z\qquot G \ar[r]^{\overline{\vphi}} & \modspc{ss}{A} }
  \end{gathered}
  \end{equation}
where the horizontal maps $\overline{i}$ and $\overline{\vphi}$ 
are closed scheme-theoretic embeddings.

In characteristic $p$, there is no Reynolds operator
and the two processes above do not commute.
Hence we cannot construct the good quotient $Z\qquot G$ as the image
$\pi_A(Z)$.
In this case, we recall that $\modspc{ss}{A}$ is a GIT quotient of $\Rep$
and $\Rep^{ss}$ is the set of GIT semistable points (\cite{Ki}).
Therefore the GIT quotient of the affine scheme $\overline{Y}$ is
a projective scheme (for the same reason as $\modspc{ss}{A}$ is \cite{Ki})
and this GIT quotient is the good quotient 
of $\overline{Y}\cap \Rep^{ss}=Z$.
We can then (uniquely) complete the diagram~\eqref{diag:goodquotbar},
because the vertical maps $\pi$ and $\pi_A$ are
categorical quotients.
We further deduce that $\overline{\vphi}$ is closed set-theoretic embedding
because $\overline{i}$ is a closed (scheme-theoretic) embedding and 
the two quotients are good.
\end{proof}

We are now in a position to complete the promised `functorial
construction' of the moduli space of semistable sheaves.

\begin{theorem}\label{thm:the-theorem}
  The scheme $\modspc{ss}{X}(\hp)$,
  constructed in Proposition~\ref{prop:goodquotient},
  is the moduli space of semistable sheaves on $X$ of Hilbert polynomial $\hp$,
  i.e. it corepresents the moduli functor $\modfun{ss}{X}(\hp)$.
  The closed points of $\modspc{ss}{X}$ correspond to 
  the S-equivalence classes of semistable sheaves.
  Furthermore, there is an open subscheme 
  $\modspc{s}{X}\subset\modspc{ss}{X}$ which corepresents
  the moduli functor $\modfun{s}{X}$ of stable sheaves and whose closed
  points correspond to the isomorphism classes of
  stable sheaves.
\end{theorem}

\begin{proof}
  To see that $\modspc{ss}{X}$ corepresents $\modfun{ss}{X}$,
  we apply Lemma~\ref{lem:loc-isom} to the local isomorphism 
  $\gtilss\colon \funpt{\Parss}/\funpt{G}\to \modfun{ss}{X}$
  of \eqref{eq:gmap-ss}
  to obtain the natural transformation
  $\psi_X$ in the following commutative diagram
\begin{equation}
\begin{gathered}\label{diag:psiX}
  \xymatrix{
    {\funpt{\Parss}}
      \ar[rd]^-{\funpt{\pi_X}}
      \ar[d]_{\gmapss} 
    & \\
    \modfun{ss}{X} 
      \ar[r]_-{\psi_X} &  
    {\funpt{\modspc{ss}{X}}}
  }
\end{gathered}
\end{equation}
and use the fact that $\pi_X$, from \eqref{diag:goodquot}, 
is a good, hence categorical, quotient.

Now, the morphism $\vphi\colon \modspc{ss}{X} \to \modspc{ss}{A}$
of \eqref{diag:goodquot} induces a bijection between the
closed points of $\modspc{ss}{X}$ and the closed points of
\[
  \pi_A\(\Parss\)\subset\modspc{ss}{A}
.\]
Hence Theorem~\ref{thm:goodquot-Kronecker} implies that the closed points
of $\modspc{ss}{X}$ correspond to the S-equivalence classes
of semistable $A$-modules that are of the form $\Hom_X(T,E)$ for semistable
sheaves $E$ of Hilbert polynomial $\hp$. 
However, we also know from
Corollary~\ref{cor:grMgrE} that $\Hom_X(T,E)$ and $\Hom_X(T,E')$ are
S-equivalent $A$-modules if and only if $E$ and $E'$ are
S-equivalent sheaves.  
Thus the closed points of $\modspc{ss}{X}$ correspond to 
the S-equivalence classes of semistable sheaves.

For the parts of the theorem concerning stable sheaves,
recall from Theorem~\ref{thm:Ess-Mss}(b) that a semistable sheaf $E$
is stable if and only if the $A$-module $\Hom_X(T,E)$ is stable.
Hence $\Pars=\Parss\cap \Rep^s$, 
where $\Rep^s$ is the open set of stable points.
In particular, all $G$-orbits in $\Pars$ are closed in $\Parss$,
because they are closed in $\Rep^{ss}$.
Therefore, we may apply Lemma~\ref{lem:goodquotient} to deduce that
$\Pars$ has a good (in fact, geometric)
quotient $\modspc{s}{X}=\Pars\qquot G=\pi_X\(\Pars\)$
which is open in $\modspc{ss}{X}=\Parss\qquot G$.
As above, $\modspc{s}{X}$ corepresents $\modfun{s}{X}$,
by Lemma~\ref{lem:loc-isom}.

Finally, the closed points of $\modspc{s}{X}$ correspond
to the isomorphism classes of stable sheaves because
$\modspc{s}{X}=\pi_X\(\Pars\)$ and 
`S-equivalence' for stable sheaves is `isomorphism'.
\end{proof}

The functorial nature of the construction can be summarised in the
following commutative diagram of natural transformations,
\begin{equation}\begin{gathered}
\label{diag:psiphi}
  \xymatrix{ 
    \modfun{ss}{X} \ar[r]^{f} \ar[d]_{\psi_X}
 &  \modfun{ss}{A} \ar[d]^{\psi_A}
 \\ {\funpt{\modspc{ss}{X}}} \ar[r]^{\funpt{\vphi}} 
 &  {\funpt{\modspc{ss}{A}}} 
  }
\end{gathered}\end{equation}
where $f$ is induced by the functor $\Hom_X(T,\?)$,
as in \eqref{diag:QRMM}, and $\psi_A$ is the corepresenting
transformation of \eqref{diag:psiA}.
If $\modspc{ss}{X}$ and $\psi_X$
are to exist, with the required properties that $\psi_X$
corepresents $\modfun{ss}{X}$ and induces a bijection between
S-equivalence classes and points of $\modspc{ss}{X}$,
then the map $\vphi$ must exist and
be an injection on points.
On the other hand, the logic of the construction is that we can effectively
show that $\psi_X$ does exist by constructing $\modspc{ss}{X}$ and $\vphi$
as in \eqref{diag:goodquot}.

\subsection{Properness of the moduli} 
\label{sub:proper}

  We already know from Proposition~\ref{prop:goodquotient} that
  $\modspc{ss}{X}$ is quasi-projective, hence to show that it is projective,
  it is sufficient to show that it is proper.
  The basic tool for this is the method of Langton 
  \cite{La} (see also \cite[\S5]{Ma02}).

\begin{theorem}\label{thm:Langton}
 Let $C$ be the spectrum of a discrete valuation ring
 and $C_0$ the generic point. 
 If $F$ is a flat family over $C_0$ of semistable sheaves on $X$, 
 then $F$ extends to a flat family of semistable sheaves over $C$.
\end{theorem}

\begin{proof}
 First (cf. the proof of \cite[Theorem 2.2.4]{HL})
 note that $F$ extends to a flat family
 over $C$, which can then be modified at the closed point,
 by \cite[Theorem 2.B.1]{HL}, to obtain a flat family of semistable sheaves.
\end{proof}

Using this we prove the following.

\begin{proposition}\label{prop:MXproper}
 The moduli space $\modspc{ss}{X}$ is proper and hence projective.
\end{proposition}

\begin{proof}
  We use the valuative criterion for properness.  
  Let $\Delta=\Spec D$ and $\Delta_0=\Spec K$, where $D$ is a
  discrete valuation ring with field of fractions $K$. 
  Given any $x_0\colon \Delta_0\to \modspc{ss}{X}$, 
  we need to show that $x_0$ extends to a map
  $x\colon \Delta\to \modspc{ss}{X}$,
  i.e. $x_0=x\circ j$,
  where $j\colon \Delta_0\hra\Delta$ is the inclusion.
  
  The first step is to lift $x_0\in \funpt{\modspc{ss}{X}}(\Delta_0)$
  along the natural transformation
  $\psi_X\colon\modfun{ss}{X}\to \funpt{\modspc{ss}{X}}$
  of \eqref{diag:psiX} to obtain a family to which we can apply
  Langton's method.

  In fact, this lift can only be achieved up to a finite cover,
  i.e. we must take a finite field extension $K'\supset K$,
  with corresponding covering map $p_0\colon \Delta_0' \to \Delta_0$,
  in order to find $y_0\in \funpt{\Parss}(\Delta'_0)$
  such that the following diagram commutes. 
\[
  \xymatrix{ 
    \Delta'_0 \ar[d]_-{p_0} \ar[r]^-{y_0}  
  & \Parss \ar[d]^-{\pi_X}
 \\ \Delta_0 \ar[r]^-{x_0}
  & \modspc{ss}{X}
  }
\]
 Then the `lift' of $x_0$ is 
 $[E_0]= \gmapss_{\Delta'_0}(y_0)\in \modfun{ss}{X}(\Delta'_0)$,
 since this family has classifying map
\begin{equation}
  \label{eq:map-E_0}
  \(\psi_X\)_{\Delta'_0}([E_0]) 
  = \(\funpt{\pi_X}\)_{\Delta'_0}(y_0) 
  = \pi_X \circ y_0
  = x_0\circ p_0. 
\end{equation}
  In other words, $E_0 = y_0^*\FF$, 
  where $\FF=i^*\MM\otimes_A T$
  is the tautological family of semistable sheaves on $\Parss$.
  
  Now, let $D'$ be a discrete valuation ring
  dominating $D$, with field of fractions $K'$
  Let $p\colon \Delta' \to \Delta$ be the corresponding
  covering map and $j'\colon \Delta'_0\hra\Delta'$ be the inclusion.
  By Theorem~\ref{thm:Langton}, $E_0$ extends over
  $\Delta'$ to a flat family $E$ of semistable sheaves on $X$, 
  i.e. $E_0=\jps E$.
  But, $E$ has classifying map
\[ 
  x' = \(\psi_X\)_{\Delta'}([E]) \colon \Delta'\to \modspc{ss}{X}
\]
  and, by the naturality of $\psi_X$,
\[
  x_0 \circ p_0 = \(\psi_X\)_{\Delta_0'}(\jps E)
  = \jps \(\psi_X\)_{\Delta'}(E) = x'\circ j'
.\]
  Since $D$ and $D'$ are discrete valuation rings and $D'$ dominates
  $D$, it follows that $D = K \cap D'$. 
  In other words, the diamond below is a push-out and therefore
  there exists a map $x\colon \Delta \to \modspc{ss}{X}$ 
  making the whole diagram commute.
\[
  \xymatrix{
    & \Delta' \ar[dr]_-{p} \ar@/^/[drrr]^-{x'}  
&&& \\ \Delta'_0 \ar[ur]^-{j'} \ar[dr]_-{p_0} 
   && \Delta \ar[rr]^-{x}
   && \modspc{ss}{X} 
  \\& \Delta_0 \ar[ur]^-{j} \ar@/_/[urrr]_-{x_0} 
   &&& 
  }
\]
  In particular, $x_0=x\circ j$, as required.
\end{proof}

\subsection{Conclusions about the embedding of moduli spaces} 
\label{sub:conc-emb}

In this subsection, we look more closely at 
the `parameter space' $\Parss\subset \Rep$ of semistable sheaves
embedded by the functor $\Hom_X(T,\?)$, 
and the induced embedding of
$\modspc{ss}{X}$ in $\modspc{ss}{A}$.

Note that, for the purposes of our construction, 
the most significant facts were that
$\Parss$ is locally closed and that $\Parss\subset \Rep^{ss}$,
which is the ``only if'' part of
Theorem~\ref{thm:Ess-Mss}(a).
Using the ``if'' part we can naturally say more.

Indeed, recall that $\Parss$ was defined as the open subset 
of semistable sheaves inside
the parameter space $\Par$ of embedded $n$-regular sheaves
and note that this is not necessarily the same as
$\Par\cap \Rep^{ss}$.
However, if we define $\Parpur$ as the open subset pure sheaves inside
$\Par$, then Theorem~\ref{thm:Ess-Mss}(a) does in fact tell us that
\begin{equation}
  \Parss = \Parpur \cap \Rep^{ss}.  
\end{equation}

In \secref{sub:modules-sheaves} we proved a stronger result,
with stronger assumptions on $n,m$.
If $\Rpur$ is the subset of $\Rep_\hp$ 
(cf. \eqref{eq:flat-strat}) consisting of
pure sheaves of Hilbert polynomial $\hp$
(not necessarily embedded), then Proposition~\ref{prop:Mss-Ess}
says that 
\begin{equation}
  \Parss = \Rpur \cap \Rep^{ss}.  
\end{equation}
Simpson's construction uses an even stronger result
\cite[Theorem 1.19]{Si},
which, in the light of Remark~\ref{rem:goodquotient-Kronecker},
is effectively that
\begin{equation}
  \Parss = \overline{\Rpur} \cap \Rep^{ss},
\end{equation}
where $\overline{\Rpur}$ is the closure in $\Rep$ 
and hence is an affine scheme. 
This is what enables Simpson to construct 
$\modspc{ss}{X}=\Parss\qquot G$ 
as the \emph{a priori} projective GIT quotient of $\overline{\Rpur}$.

By contrast, in this paper, we used Langton's method
to  discover \emph{a posteriori} that 
$\Parss\qquot G$ is proper.
In the notation of Proposition~\ref{prop:goodquotient},
i.e. $Y=\Parss$ and $Z=\overline{Y} \cap \Rep^{ss}$, 
we may then deduce that $Y\qquot G$ is closed in 
the (separated) GIT quotient $Z\qquot G$.
Therefore, as $Y\qquot G$ is also dense, it is equal to $Z\qquot G$. 
Since $Y=\pi^{-1}(\pi(Y))$ 
(cf. Lemma~\ref{lem:goodquotient}),
this means that $Y=Z$, i.e.
\begin{equation}
\label{eq:Y=Z}
  \Parss = \overline{\Parss} \cap \Rep^{ss}.
\end{equation}
Thus, we now could say, with weaker assumptions on $n,m$ than Simpson
needs, 
that $\modspc{ss}{X}$ is the GIT quotient of 
the affine scheme $\overline{\Parss}\subset \Rep$.

In conclusion, we have

\begin{proposition}
  \label{prop:embedding-final}
 The map $\vphi\colon\modspc{ss}{X}(\hp)\to \modspc{ss}{A}\(\hp(n),\hp(m)\)$
 in \eqref{diag:goodquot} is a closed set-theoretic embedding of 
 projective schemes. 
 This embedding is
 scheme-theoretic in characteristic zero, while
 in characteristic $p$ it is scheme-theoretic on the stable locus.
 
 Furthermore, the GIT construction yields an
 ample line bundle $L$ on $\modspc{ss}{A}$ such that $\vphi^*L$
 is an ample line bundle on $\modspc{ss}{X}$.
\end{proposition}

\begin{proof}
First recall that $\modspc{ss}{A}$ is projective by construction,
while $\modspc{ss}{X}$ is projective by Proposition~\ref{prop:MXproper}.
From above, \eqref{eq:Y=Z} means that $\vphi$ coincides with 
the map $\overline{\vphi}$ in \eqref{diag:goodquotbar} 
and hence is a closed embedding 
(scheme-theoretic in characteristic zero, 
but only set-theoretic in characteristic $p$).
To see that $\vphi$ is a scheme-theoretic embedding on the stable locus
even in characteristic $p$, consider the restriction of the diagram
\eqref{diag:goodquot} to the stable loci.
\begin{equation}\begin{gathered}
\label{diag:stable-goodquot}
  \xymatrix{ 
    \Pars \ar[r]^{i} \ar[d]_{\pi_X}
 &  \Rep^s \ar[d]^{\pi_A}
 \\ \modspc{s}{X} \ar[r]^{\vphi} 
 &  \modspc{s}{A} 
  }
\end{gathered}\end{equation}
In \eqref{diag:goodquot} we know that $i$ is a closed (scheme-theoretic)
embedding by \eqref{eq:Y=Z}.
Hence, in \eqref{diag:stable-goodquot} we also know that 
$i$ is a closed embedding, because $\Pars=\Parss\cap \Rep^s$ 
by Theorem~\ref{thm:Ess-Mss}(b).
On the other hand, we know that a stable module (or a stable sheaf) 
is `simple', in the sense that its endomorphism algebra is just $\kk$.
Thus $G$ acts freely on $\Rep^s$ and so 
$\pi_A$ is a principal $G$-bundle over $\modspc{s}{A}$.
But a closed $G$-invariant subscheme of a principal bundle
is a principal bundle over a closed subscheme of the base and thus
the restricted $\vphi$ in \eqref{diag:stable-goodquot} is a 
scheme-theoretic closed embedding.

Finally, note that $\modspc{ss}{A}$ is constructed as a GIT quotient of $\Rep$,
with respect to a $G$-linearised (trivial) 
line bundle $\cL$, for which some power $\cL^N$ descends from $\Rep^{ss}$
to an ample line $L$ on $\modspc{ss}{A}$.
On the other hand, we now know that $\modspc{ss}{X}$ can be constructed
as a GIT quotient of $\overline{\Parss}$ with respect to the restriction
of $\cL$ and, by the diagram \eqref{diag:goodquot},
the restriction of $\cL^N$ to $\Parss$ descends to $\vphi^* L$,
which is therefore ample.
\end{proof}

In Section~\ref{sec:cohochar} we will look more closely at
the line bundle $\cL$ (and its $G$-invariant sections) and 
we will see that, in fact, all powers of $\cL$ descend.

One could also deduce that $\vphi^* L$ is ample
directly from the ampleness of $L$, because $\vphi$ is a finite map
(cf.  \cite[Ch III, Ex 5.7]{Ha}).

\subsection{Uniform/universal properties of moduli spaces} 
\label{sub:uni-properties}

The moduli spaces of stable and semistable sheaves actually enjoy 
a stronger property than corepresenting the moduli functors.

\begin{definition}\label{def:univ/unif}
  Let $\MMM\colon \Sch^\circ\to \Set$ be a functor, $\cM$ a scheme
  and $\psi\colon \MMM\to\funpt{\cM}$ a natural transformation.
  We say that $\cM$ \emph{universally
  corepresents} $\MMM$ if for each morphism of schemes $\cN\to\cM$,
  the fibre product functor $\funpt{\cN} \times_{\funpt{\cM}}
  \MMM$ is corepresented by the canonical projection $\funpt{\cN}
  \times_{\funpt{\cM}} \MMM \to \funpt{\cN}$. If this holds
  only for \emph{flat} morphisms $\cN\to\cM$, then we say that $\cM$
   \emph{uniformly corepresents} $\MMM$.
\end{definition}

This enhancement of Definition~\ref{def:corep}
is a direct generalisation of a similar
enhancement of `categorical quotient' in \cite[Definition~0.7]{Mu1}.
Indeed, if an algebraic group $G$ acts on a scheme $Z$, 
then a $G$-invariant morphism $Z\to Z\qquot G$ is a
uniform/universal categorical quotient if and only if the induced natural
transformation $\funpt{Z}/ \funpt{G} \to \funpt{Z\qquot G}$
uniformly/universally corepresents $\funpt{Z}/ \funpt{G}$.
In the situation of Lemma~\ref{lem:loc-isom},
we can conclude that $\psi_1$ uniformly/universally corepresents $\AAA_1$ if
and only if $\psi_2$ uniformly/universally corepresents $\AAA_2$.

We can enhance Lemma~\ref{lem:goodquotient} by adding
``uniform'' (or ``universal'' in characteristic zero)
to ``good quotient'', because 
\[ 
  Y=\pi^{-1}(Y\qquot G)=Y\qquot G\times_{Z\qquot G} Z
.\]
We can similarly enhance Proposition~\ref{prop:goodquotient}.
In the characteristic zero case, this is because $\pi_A$ is 
a universal good quotient, because it is a GIT quotient.
This property is then automatically inherited by $\pi$,
because the diagram \eqref{diag:goodquotbar} is a pull-back.
In the characteristic $p$ case, we must be more direct:
$\pi$ is a uniform good quotient, 
because it is a GIT quotient.

Thus we can enhance Theorem~\ref{thm:the-theorem} by adding
``uniformly'' (or ``universally'' in characteristic zero)
to ``corepresents''.

\subsection{Relative moduli spaces} 
\label{sub:relative-moduli}

The functorial method also enables us to construct relative moduli
spaces for a projective morphism of schemes $\rho\colon X\to Y$, with
a relatively very ample invertible sheaf $\OX(1)$.

If $\Sch/Y$ is the category of schemes over $Y$, then we must work
with the relative moduli functors
\[
 \modfun{ss}{X}, \modfun{s}{X} 
 \colon (\Sch/Y)\op\to \Set
\]
defined as in~\secref{sub:construction-moduli}, where now a `flat
family over $S$ of sheaves on $X$' is an $S$-flat sheaf on $X\times_Y
S$.

Let $T = \OX(-n) \oplus \OX(-m)$ and, for $H=\rho_*(\OX(m-n))$, let
\[
 A=\begin{pmatrix} \cO_Y & H \\ 0 & \cO_Y\end{pmatrix},
\]
an $\cO_Y$-algebra which naturally acts on $T$.
In other words, there is an $\cO_Y$-algebra morphism
$A\to\rho_*\cEnd_X(T)$, defined using the structure map $\cO_Y\to
\rho_*\cO_X$, the projection operators and the obvious $H$-action. 
Thus, we have a functor from sheaves of $\cO_X$-modules on $X$ to
sheaves of $A$-modules on $Y$, written simply $\Hom_X(T,\?) \defeq
\rho_*\cHom_X(T,\?)$. 
This functor also has a left adjoint, written $\?\otimes_A T$.

Now 
a sheaf $E$ on $X$
is \emph{$n$-regular} (relative to $\rho$) if
$R^i\rho_*(E(n-i))=0$ for all $i>0$.
Then (cf.~\cite[Ch V, Prop 2.2]{FL})  
there is a relative version of Lemma~\ref{lem:regular}
and hence relative versions of  Lemma~\ref{lem:regker}
and Theorem~\ref{thm:eval=isom}. 

An additional technical point in the relative case, is that
we should (and can) choose $m-n$ large enough that any
base change of $\OX(m-n)$ is regular and 
further the formation of $\rho_*(\OX(m-n))$ commutes with base
change (see~\cite[Lemma~1.9]{Si} and the remark following it).
Under this additional assumption, there are relative versions of
Proposition~\ref{prop:fullyfaithful} and hence
Proposition~\ref{prop:nreg-image}, where a flat family of $A$-modules
over a relative scheme $\sigma\colon S\to Y$ is a locally free sheaf
of right modules over $\sigma^*(A)$.
The relative version of condition (C:3) 
in~\secref{sub:long-list} should be enhanced by adding this assumption. 

Now we have relative moduli functors
\[
 \modfun{ss}{A}\subset\modfun{}{A}
 \colon (\Sch/Y)\op\to \Set
\]
defined as in~\secref{sub:mod-fun-Kron-mod}, but using flat families
of $A$-modules over relative schemes.
%
To construct a relative moduli space of Kronecker modules (as in
Theorem~\ref{thm:goodquot-Kronecker}), choose free (or even locally
free) sheaves $V$ and $W$ over $Y$ of ranks $\hp(n)$ and $\hp(m)$,
respectively.
There is a scheme $\Rep$ over $Y$, which parametrises representations
of $A$ on $V\oplus W$, or equivalently morphisms $H\to \cHom_Y(V,W)$.
In other words, $\Rep$ represents the functor $(\Sch/Y)\op\to \Set$
which assigns to any $\sigma\colon S\to Y$ the set $\Hom_S(\sigma^*
H, \sigma^*\cHom_Y(V,W))$. Note that the existence of this scheme
depends on the fact that $\cHom_Y(V,W)$ is locally
free~\cite[Theorem~3.5]{Ni}.

The group (scheme over $Y$) $G=\GL(V)\times\GL(W)/\Delta$ acts
naturally on $\Rep$.  If $\Rep^{ss}\subset \Rep$ is the open set of
points corresponding to semistable modules over the fibres of $A$,
then the quotient functors $\funpt{\Rep^{ss}}/\funpt{G} \subset
\funpt{\Rep}/\funpt{G}$ are locally isomorphic (over $Y$) to
$\modfun{ss}{A} \subset \modfun{}{A}$, as in~\secref{sub:mod-fun-Kron-mod}. 
Now GIT constructs a scheme $\modspc{ss}{A}$, projective over $Y$,
and a (relative) good quotient $\Rep^{ss}\to\modspc{ss}{A}$ for the
$G$-action. 
Thus $\modspc{ss}{A}$ corepresents $\modfun{ss}{A}$.

Applying the relative version of Proposition~\ref{prop:nreg-image} to
the tautological family $\MM$ of right $A$-modules on $\Rep$, we
obtain a $Y$-scheme $\Par = \Rreg\subset \Rep$, which satisfies a
relative version of Theorem~\ref{thm:reg-loc-isomX}. 
Choosing $n$ large enough that condition (C:1) of
\secref{sub:long-list} is satisfied for the fibre of $\rho$ over
every point of $Y$, we can define $\Pars\subset \Parss \subset \Par$
and obtain a relative version of Proposition~\ref{prop:sss-loc-isomX}. 
Choosing $n$ and $m$ large enough that conditions (C:2), (C:4)
and (C:5) of \secref{sub:long-list} are also satisfied for the fibre
of $\rho$ over every point of $Y$, Theorem~\ref{thm:Ess-Mss} and
Corollary~\ref{cor:grMgrE} hold fibrewise.
As the formation of $\rho_*(\OX(m-n))$ commutes with base change, the 
fibre of $H$ over a point of $Y$ is the space of sections of $\OX(m-n)$
restricted to the corresponding fibre of $\rho$. 
Thus, Theorem~\ref{thm:Ess-Mss} and 
Corollary~\ref{cor:grMgrE} apply to the fibres of $\Hom_X(T,\FF)$ over the 
points of $Y$, with $\FF$ as in Proposition~\ref{prop:sss-loc-isomX}.
Then we see that $\Parss$ is locally closed in $\Rep^{ss}$
and we construct $\modspc{ss}{X}$ as a relative good quotient
$\Parss\qquot G$, as in Proposition~\ref{prop:goodquotient}.

Hence we obtain a relative version of Theorem~\ref{thm:the-theorem}
in which $\modspc{ss}{X}$ is a scheme over $Y$ corepresenting the
relative moduli functor $\modfun{ss}{X}$ and the closed points of
$\modspc{ss}{X}$ correspond to the S-equivalence classes of
semistable sheaves over the fibres of $\rho$.  There is a similar
modification for $\modspc{s}{X}$.
In addition, the proof of Proposition~\ref{prop:MXproper} applies in
the relative case, so $\modspc{ss}{X}$ is proper, and hence
projective, over $Y$.

Finally, note that this enhancement can be combined with the
enhancement in \secref{sub:uni-properties} by adding ``uniformly''
(or ``universally'' in characteristic zero) to ``corepresents''.

\section{Determinantal theta functions}
\label{sec:cohochar}

In this section, we interpret the main results of the paper in terms
of determinantal theta functions on the moduli space $\modspc{ss}{X}$
of sheaves,
using analogous results already proved by 
Schofield \& Van den Bergh \cite{SV} and Derksen \& Weyman \cite{DW}
for the moduli space $\modspc{ss}{A}$ of Kronecker modules,
or more generally, for moduli spaces of representations of quivers.
A key ingredient is the adjunction that has been central to this paper,
between $\Phi=\Hom_X(T,\?)$ and $\Phiadj=\?\otimes_A T$.

In \secref{sub:det-char-ss}, we give a new characterisation of 
semistable sheaves, amongst $n$-regular pure sheaves $E$,
as those which invert certain maps of vector bundles.
The condition is implicitly equivalent to the semistability of
$\Hom_X(T,E)$, but without any explicit reference to Kronecker modules.
This means that semistable sheaves are characterised by the non-vanishing of
the corresponding determinantal theta functions, which we describe in more
formal detail in \secref{sub:theta-on-moduli}, showing in particular
how to interpret them as sections of line bundles on the moduli space.
Using stronger results of \cite{SV,DW}, we show in
\secref{sub:separation-prop}, that determinantal theta functions
can actually be used to give a projective embedding
the moduli space $\modspc{ss}{X}$, modulo the technical problems
with semistable points in characteristic $p$ that we have already
encountered.
Finally, in \secref{sub:char-curve}, we explain
how the results of \secref{sub:separation-prop} improve what
was known even in the case when $X$ is a smooth curve.

\subsection{A determinantal characterisation of semistability} 
\label{sub:det-char-ss}

To see how semistable $A$-modules can be characterised,
note that, as a right $A$-module $A=P_0\oplus P_1$, where
$P_0=e_0A$ and $P_1=e_1A$
are the indecomposable projective modules.
If $M$ is an $A$-module, then the corresponding Kronecker
module $\alpha\colon V\otimes H\to W$ is given by the
composition map, after noting that
\[
   V=\Hom_A(P_0,M),
\quad W=\Hom_A(P_1,M),
\quad H=\Hom_A(P_1,P_0).
\]
Now, a corollary of the results of \cite{SV,DW}
can be formulated as saying that semistable $A$-modules
are precisely those which invert certain maps between projective modules.

\begin{theorem}
\label{thm:SVDW}
An $A$-module $M=V\oplus W$ is semistable if and only if there is a map
\begin{equation}
\label{eq:gamma}
 \gamma\colon U_1\otimes P_1 \to U_0\otimes P_0,
\end{equation}
where $U_i$ are (non-zero) vector spaces,
such that the induced linear map 
\[
  \Hom_A(\gamma,M)\colon \Hom(U_0,V)\to \Hom(U_1,W)
\]
is invertible,
i.e. $\theta_\gamma(M)\defeq\det\Hom_A(\gamma,M)\neq0$.
\end{theorem}

\begin{proof}
First note that, the fact that $\Hom_A(\gamma,M)$ may be invertible
requires in particular that
\begin{equation}
\label{eq:UVW}
  \dim U_0 \dim V = \dim U_1 \dim W.  
\end{equation}
When this holds, $\theta_\gamma$ is a $G$-equivariant polynomial function on
the representation space $\Rep$ of \eqref{eq:rep-space},
with values in the one-dimensional $G$-vector space
\[
  (\det V)^{-\dim U_0} \otimes (\det W)^{\dim U_1}.
\]
Identifying this space with $\kk$, we can consider $\theta_\gamma$
as a semi-invariant with weight
\begin{equation}
\label{eq:chiU}
  \chi_U \colon G \to \kk^\times 
  \colon (g_0,g_1)
  \mapsto \(\det g_0\)^{-\dim U_0} \(\det g_1\)^{\dim U_1},
\end{equation} 
that is, $\theta_\gamma (g\cdot \alpha)=\chi_U(g)\theta_\gamma (\alpha)$ 
for all $\alpha\in \Rep$ and $g\in G$. 
By \cite{Ki}, GIT-semistable points $\alpha\in \Rep$ with respect to the
character $\chi_U$ correspond one-one with
semistable $A$-modules $M$, because the condition from \cite{Ki}
on submodules $M'=V'\oplus W'$ that
\[
  \dim U_1 \dim W' - \dim U_0 \dim V' \geq 0,
\]
is equivalent to the condition \eqref{eq:55} by \eqref{eq:UVW}.
In other words, $M$ (or $\alpha$) is semistable if and only if there is some
semi-invariant $\theta$ with $\theta(\alpha)\neq 0$,
where $\theta$ has weight $\chi_U$ for some $U_0,U_1$ satisfying
\eqref{eq:UVW}.
But, by \cite[Theorem~2.3]{SV} or  \cite[Theorem~1]{DW},
the space of semi-invariants of weight $\chi_U$ is spanned by the `determinantal
semi-invariants' of the form $\theta_\gamma$ and so the result follows.

More precisely, in the notation of \cite[\S 3]{SV},
we have $\theta_\gamma=P_\phi$, where $\gamma=\hat\phi$,
and, in the notation of \cite{DW},
we have $\theta_\gamma=c^N$, where $\gamma$ is a projective resolution of $N$.

It is an extra observation, in this case, that 
$\gamma$ may be chosen with the particular domain and codomain of
\eqref{eq:gamma}.
From the perspective of \cite{DW}, this occurs because $N$ 
must also be semistable and hence, in particular, saturated. 
From the perspective of \cite{SV}, it is an elementary computation
that inverting a general map between projective $A$-modules
is equivalent to inverting a map of this form.
\end{proof}

From this theorem for $A$-modules, we may derive a similar
determinantal characterisation of semistability for sheaves on $X$.

\begin{theorem}
\label{thm:cohochar-ssbundles}
For a fixed polynomial $\hp$, suppose that $n,m$ satisfy
the conditions (C:1)-(C:5) of \secref{sub:long-list}.

Let $E$ be an $n$-regular pure sheaf of Hilbert polynomial $\hp$.
Then $E$ is semistable if and only if
there is a map
\begin{equation}
\label{eq:delta}
  \delta\colon U_1\otimes\OX(-m) \to U_0\otimes\OX(-n)
\end{equation}
where $U_i$ are (non-zero) vector spaces,
such that the induced linear map 
\[
  \Hom_X(\delta,E)\colon \Hom(U_0,H^0(E(n)))\to \Hom(U_1,H^0(E(m)))
\]
is invertible,
i.e. $\theta_\delta(E)\defeq\det\Hom_X(\delta,E)\neq0$.
\end{theorem}

\begin{proof}
As $\Phiadj(A)=T$, so $\Phiadj(P_0)=\OX(-n)$ and 
$\Phiadj(P_1)=\OX(-m)$.
Thus $\Phiadj$ gives a bijection between
the maps $\delta$ in \eqref{eq:delta} and the maps $\gamma$ in \eqref{eq:gamma}.
Furthermore, the adjunction between $\Phi$ and $\Phiadj$ implies that
\begin{equation}
\label{eq:adj-on-hom}
  \Hom_X(\Phiadj(\gamma),E)=\Hom_A(\gamma, \Phi(E)).
\end{equation}
Thus, the existence of $\delta$ with $\theta_\delta(E)\neq0$
is equivalent to the semistability of $\Phi(E)$,
by Theorem~\ref{thm:SVDW}, which is equivalent to the semistability of $E$,
by Theorem~\ref{thm:Ess-Mss}(a).
\end{proof}

\begin{remark}
\label{rem:elem-proof}
The ``if'' part of Theorem~\ref{thm:cohochar-ssbundles}
has a more direct and elementary proof.
Note first that, if such a $\delta$ exists, 
then, as the domain and codomain of $\Hom_X(\delta,E)$ 
must certainly have the same dimension, we know that
\begin{equation}
  \label{eq:euler=0}
  \hp(n)\dim U_0 = \hp(m) \dim U_1,
\end{equation}
because $E$ is $n$-regular.
Now, if $E$ were not semistable,
then by Proposition~\ref{prop:new-est-E},
there would exist a subsheaf $E'\subset E$ with 
\[
  h^0(E'(n)) \hp(m) > \hp(n) h^0(E'(m))
\]
and thus 
\[
  h^0(E'(n)) \dim U_0 > h^0(E'(m)) \dim U_1. 
\]
In other words, if we write $K_0=U_0\otimes\OX(-n)$ 
and $K_1=U_1\otimes\OX(-m)$, then 
\[
  \dim\Hom_X(K_0,E') > \dim\Hom_X(K_1,E').
\]
Hence the map 
\[
  \Hom_X(\delta,E')\colon \Hom_X(K_0,E')\to\Hom_X(K_1,E')
\]
has non-zero kernel, i.e. there is a non-zero map 
$\phi\colon K_1\to E'$, with $\phi\circ \delta=0$.
But, as $E'\subset E$, this also shows that the kernel of $\Hom_X(\delta,E)$
is non-zero, contradicting the assumption.
Thus $E$ must be semistable.
\end{remark}

\begin{remark}
\label{rem:dercat}
One may interpret Theorem~\ref{thm:cohochar-ssbundles} in terms
of derived categories, using the derived adjunction
\begin{equation}
\label{eq:deriv-adj}
  \R\Hom_X(\L\Phiadj(N),E)
 \cong
 \R\Hom_{A}(N,\R\Phi(E)).
\end{equation}
As $E$ is $n$-regular, $\R\Phi(E)=\Phi(E)$,
while, by \cite[Theorem~1]{DW} and \cite[Corollary~1.1]{SV},
the semistability of $\Phi(E)$ is equivalent to the existence of
an $A$-module $N$, which is `perpendicular' to $\Phi(E)$ in the sense
that the right-hand side of \eqref{eq:deriv-adj} vanishes.
On the other hand, the complex $K_\bullet$ that represents
$\L\Phiadj(N)$ is obtained by taking a projective
resolution of $N$ given by a map $\gamma$ as in \eqref{eq:gamma}
and applying $\Phiadj$.
In other words, we obtain the map $\delta$ as in \eqref{eq:delta},
interpreted as a 2-step complex
\[
  K_\bullet = K_1\lra{\delta} K_0.
\]
Observe that, again as $E$ is $n$-regular, the perpendicularity
 condition
\[
  \R\Hom_X(K_\bullet,E)=0
\] 
is precisely the condition 
$\theta_\delta(E)\neq0$ of Theorem~\ref{thm:cohochar-ssbundles}.
\end{remark}

\subsection{Theta functions on moduli spaces} 
\label{sub:theta-on-moduli}

We now explain how the `functions'
$\theta_\gamma$ and $\theta_\delta$ of \secref{sub:det-char-ss}
can be properly interpreted as sections of line bundles
on the moduli spaces $\modspc{ss}{A}$ and $\modspc{ss}{X}$
with certain universal properties, describing first
in some detail the case of $A$-modules.

Consider a flat family $M=V\oplus W$ over $S$ of
$A$-modules with dimension vector $(a,b)$.
Then for any $\gamma\colon U_1\otimes P_1 \to U_0\otimes P_0$,
with
\begin{equation}\label{eq:abU}
  a \dim U_0=b\dim U_1,
\end{equation}
we may define a line bundle over $S$
\[
 \lambda_U(M) \defeq\(\det\Hom(U_0,V)\)^{-1}\otimes\det\Hom(U_1,W)
\]
and a global section
\[
 \theta_\gamma(M)\defeq\det\Hom_A(\gamma,M)\in H^0(S,\lambda_U(M)).
\]
Roughly, we have a natural assignment, to each module $M$, 
of a one-dimensional vector space together with a vector in it.
More formally, in the sense of \cite[\S 3.1]{DN},
we have a line bundle $\lambda_U$, together with a global
section $\theta_\gamma$, on the moduli functor $\modfun{}{A}(a,b)$
of~\eqref{eq:mod-funA}.
This means that, given another family $M'$ over $S'$
such that $M'\cong \sigma^*M$ for $\sigma\colon S'\to S$,
there is an isomorphism $\sigma^*\lambda_U(M)\cong\lambda_U(M')$
which identifies $\sigma^*\theta_\gamma(M)$ with $\theta_\gamma(M')$.
Furthermore, these identifications are functorial in $\sigma$.
 
What we can show is that the restriction of this formal line bundle and
section to the moduli functor $\modfun{ss}{A}$ of~\eqref{eq:sub-mod-funA}
descends to a genuine line bundle and section on the moduli space
$\modspc{ss}{A}$, in the following sense. 

\begin{proposition}
\label{prop:descentA}
  There is a unique line bundle $\lambda_U(a,b)$ on the moduli space 
 $\modspc{ss}{A}(a,b)$ and a global section $\theta_\gamma(a,b)$
 such that, for any family $M$ over $S$ of semistable $A$-modules 
 of dimension vector $(a,b)$, we have
\[
 \lambda_U(M)\cong \psi_M^*\lambda_U(a,b),
\qquad
 \theta_\gamma(M)=\psi_M^*\theta_\gamma(a,b),
\]
where $\psi_M\defeq \(\psi_A\)_S([M])\colon S\to\modspc{ss}{A}(a,b)$ is the
classifying map coming from \eqref{diag:psiA}. 
\end{proposition}

\begin{proof}
  Let $\MM$ be the tautological family of $A$-modules on $\Rep^{ss}(a,b)$,
as in \secref{sub:mod-fun-Kron-mod}.
Then \eqref{eq:abU} implies that $\lambda_U(\MM)$ is a 
$G$-linearised line bundle.
Thus, using Kempf's descent condition \cite[Theorem~4.2.15]{HL},
we see that $\lambda_U(\MM)$ descends to a (unique) line bundle 
$\lambda_U(a,b)$ on the good quotient
$\Rep^{ss}(a,b)\qquot G=\modspc{ss}{A}(a,b)$ if and only if
 for each point $\alpha\in \Rep^{ss}$ in a
closed orbit, the isotropy group of $\alpha$ acts trivially on the fibre 
over $\alpha$.
Note that, as $\pi_A=\psi_{\MM}$, if it exists, then $\lambda_U(a,b)$
must be the descent of $\lambda_U(\MM)$.

By \cite[Proposition~3.2]{Ki}, a point
$\alpha$ of $\Rep^{ss}$ is in a closed orbit if and only if the module
$M=\MM_\alpha$ is `polystable', that is,
\[
  M\cong \bigoplus_{i=1}^k K_i\otimes M_i,
\]
where $M_i$ are non-isomorphic stable modules of
dimension vector $(a_i,b_i)$ with the same slope as $M$
and $K_i$ are multiplicity vector spaces. 
Since stable modules are simple, the isotropy group of
$\alpha$ for the action of $\GL(V)\times \GL(W)$ is
isomorphic to
\[
  \Aut M\cong \prod_{i=1}^k \GL(K_i).
\]
By standard properties for determinants of sums and tensor products,
\[
  \lambda_U(M) \cong \bigotimes_{i=1}^k \lambda_U (K_i\otimes M_i)
  \cong 
  \bigotimes_{i=1}^k (\det K_i)^{\nu_i}\otimes
  \lambda_U(M_i)^{\dim K_i}, 
\]
as linear representations of the isotropy group, where
\[
  \nu_i \defeq b_i\dim U_1 - a_i\dim U_0 =0,
\]
by \eqref{eq:abU}, because $a_i/b_i=a/b$.
Thus, the isotropy group acts trivially, as required.

Furthermore, $\theta_\gamma(\MM)$ is a $G$-invariant section of 
$\lambda_U(\MM)$ and so it descends to a section $\theta_\gamma(a,b)$
of $\lambda_U(a,b)$, because the descent means that 
$\lambda_U(a,b)$ is the $G$-invariant push-forward of $\lambda_U(\MM)$.

The universal properties of $\lambda_U(a,b)$ and $\theta_\gamma(a,b)$
follow now from a careful analysis of the local isomorphism
$\funpt{\Rep^{ss}}/\funpt{G} \to \modfun{ss}{A}$ in~\eqref{eq:loc-isom-ssA}
along similar lines to \cite[\S 3.2]{DN}. 
\end{proof}

Note that the notation $\lambda_U$ emphasizes the fact that 
the line bundle depends on the pair of vector spaces
$U_0,U_1$, but not on the map $\gamma$.
Note also that 
\begin{equation}
  \label{eq:lamUU'}
  \lambda_{U\oplus U'}=\lambda_U \otimes \lambda_{U'},
\qquad
  \theta_{\gamma\oplus\gamma'}=\theta_{\gamma}\theta_{\gamma'}.
\end{equation}
This applies equally to the line bundles $\lambda_U(M)$
on families and the line bundles $\lambda_U(a,b)$ on the moduli space
$\modspc{ss}{A}(a,b)$, because the later descend from a specific
case of the former and pull-back commutes with tensor product. 

Indeed, up to isomorphism $\lambda_U$ depends only on $\dim U_0$ and 
$\dim U_1$ and so, with the constraint \eqref{eq:abU}, 
all possible $\lambda_U$ are isomorphic 
to positive powers of a single $\lambda_U$
with $\dim U_0$ and $\dim U_1$ coprime.

\begin{proposition}
\label{prop:lamample}
  The line bundle $\lambda_U(a,b)$ on $\modspc{ss}{A}(a,b)$ is
  ample. 
  Furthermore, its space of global sections is canonically isomorphic to
  the space of semi-invariants on $\Rep$ with the weight
  $\chi_U$ of \eqref{eq:chiU}.
\end{proposition}

\begin{proof}
Let $\MM$ be the tautological family of $A$-modules on the whole of $\Rep$.
As shown in the proof of Theorem~\ref{thm:SVDW}, $\lambda_U(\MM)$
is the $G$-linearised line bundle used in~\cite{Ki}
to construct $\modspc{ss}{A}$ as a GIT quotient,
i.e. $\lambda_U(\MM)$ is the line bundle $\cL$ 
in the proof of Proposition~\ref{prop:embedding-final}.
Thus, as the restriction of $\lambda_U(\MM)$ to $\Rep^{ss}$ does descend
to the quotient, by Proposition~\ref{prop:descentA}, 
the descended line bundle $\lambda_U(a,b)$ is ample 
and can be taken to be the line bundle $L$ 
in the statement of Proposition~\ref{prop:embedding-final}. 

As also shown in the proof of Theorem~\ref{thm:SVDW},
the semi-invariants on $\Rep$ with the weight $\chi_U$
are identified with the $G$-invariant sections of $\lambda_U(\MM)$.
As $\Rep$ is normal (cf. \cite[Theorem 4.1(ii)]{Se3}), 
we have canonical isomorphisms
\[
  H^0\(\Rep,\lambda_U(\MM)\)^G = 
  H^0\(\Rep^{ss},\lambda_U(\MM)\)^G = 
  H^0\(\modspc{ss}{A}(a,b),\lambda_U(a,b)\).
\] 
\end{proof}

Note that this last canonical isomorphism identifies the sections
$\theta_\gamma(a,b)$ of Proposition~\ref{prop:descentA}
with the corresponding determinantal semi-invariants
$\theta_\gamma(\MM)$.

Most of the above carries over similarly to the case of sheaves.
Given a map $\delta\colon U_1\otimes\OX(-m) \to U_0\otimes\OX(-n)$,
where the vector spaces $U_0,U_1$ satisfy
\begin{equation}
\label{eq:UP}
  \hp(n)\dim U_0 = \hp(m) \dim U_1,
\end{equation}
we obtain a line bundle $\lambda_U$ with a section $\theta_\delta$
on the moduli functor $\modfun{reg}{X}(\hp)$ of $n$-regular sheaves
with Hilbert polynomial $\hp$.
This is defined in the analogous way, i.e. given a family $E$ over $S$
of $n$-regular sheaves with Hilbert polynomial $\hp$, 
we have a line bundle over $S$
\[
 \lambda_U(E) \defeq \(\det\Hom_X(U_0,E(n))\)^{-1}\otimes\det\Hom_X(U_1,E(m)) 
\]
and a global section
\[
  \theta_\delta(E)\defeq\det\Hom_X(\delta,E)\in H^0(S,\lambda_U(E)),
\]
with the appropriate functorial properties.

Note that every such $\delta$ is of the form $\Phiadj(\gamma)$,
for $\gamma=\Phi(\delta)$ and hence 
the adjunction between $\Phi$ and $\Phiadj$ yields the identification
\[
  \Hom_X(\delta,E)=\Hom_A(\gamma, M),
\]
for $M=\Phi(E)$.
Hence we naturally have
\begin{equation}
\label{eq:lambda-theta}
  \lambda_U(E)\cong\lambda_U(M),
\qquad
  \theta_\delta(E)=\theta_\gamma(M).
\end{equation}

This essentially tells us how the theta functions restrict under the embedding
of Proposition~\ref{prop:embedding-final},
\[
  \vphi\colon \modspc{ss}{X}(\hp) \to \modspc{ss}{A}(a,b),
\]
where $(a,b)=(\hp(n),\hp(m))$ and $n,m$ satisfy the conditions
(C:1)-(C:5) of \secref{sub:long-list}.
More precisely, we can use Propositions~\ref{prop:descentA}
and~\ref{prop:lamample} to obtain an analogous result for
$\modspc{ss}{X}$.

\begin{proposition}
 \label{prop:descentX}
 There is an ample line bundle
 $\lambda_U(\hp)=\vphi^*\lambda_U(a,b)$ 
 on the moduli space $\modspc{ss}{X}(\hp)$ and a global section 
 $\theta_\delta(\hp)=\vphi^*\theta_\gamma(a,b)$,
 for $\gamma=\Phi(\delta)$,
 such that, for any family $E$ over $S$ of semistable sheaves with Hilbert
 polynomial $\hp$, we have
\[
 \lambda_U(E)\cong \psi_E^*\lambda_U(\hp),
\qquad
 \theta_\delta(E)=\psi_E^*\theta_\delta(\hp),
\]
where $\psi_E\defeq \(\psi_X\)_S([E])\colon S\to\modspc{ss}{X}(\hp)$ is the
classifying map coming from \eqref{diag:psiX}. 
\end{proposition}

\begin{proof}
 Firstly, $\lambda_U(\hp)$ is ample, 
 because $\lambda_U(a,b)$ is the ample line bundle 
 $L$ of Proposition~\ref{prop:embedding-final}:
 see the proof of Proposition~\ref{prop:lamample}.

 Recall further that $\modspc{ss}{X}(\hp)$ was constructed as
 a good quotient of a subscheme $\Parss\subset \Rep^{ss}$
 as in \eqref{diag:goodquot} which carries a tautological family 
 $\FF$ of semistable sheaves, with $\Phi(\FF)\cong i^*\MM$.
 Using \eqref{eq:lambda-theta} and \eqref{diag:goodquot}, 
 we see that
 \begin{align*}
   \lambda_U(\FF) &= i^*\lambda_U(\MM)=\pi_X^*\vphi^*\lambda_U(a,b)\\
   \theta_\delta(\FF) &=i^*\theta_\gamma(\MM)=\pi_X^*\vphi^*\theta_\gamma(a,b)
 \end{align*}
 that is, $\lambda_U(\FF)$ and $\theta_\delta(\FF)$ descend
 to $\lambda_U(\hp)$ and $\theta_\delta(\hp)$ as defined in the proposition.
 The universal properties now follow,
 as in the proof of Proposition~\ref{prop:descentA}, 
 or by direct argument from \eqref{diag:psiphi}.
\end{proof}

\subsection{The separation property} 
\label{sub:separation-prop}

We now use the full force of the results of \cite{DW,SV}
to obtain stronger
results about the determinantal theta functions of sheaves
$\theta_\delta$.
The point is that the
determinantal theta functions of modules $\theta_\gamma$
don't just detect semistable $A$-modules,
they actually span the ring of semi-invariants on $\Rep$, which means in particular
that they furnish a full set of homogeneous coordinates on the moduli space
$\modspc{ss}{A}$.

\begin{theorem}
\label{thm:sep-Amod}
  For any dimension vector $(a,b)$, we can find
  vector spaces $U_0,U_1$ satisfying \eqref{eq:abU}
  and finitely many maps
\[
  \gamma_0,\dots,\gamma_N \colon U_1\otimes P_1 \to U_0\otimes P_0
\]
  such that the map 
\begin{equation}
  \label{eq:Thetagam}
  \Theta_\gamma \colon \modspc{ss}{A}(a,b) \to \PP^N\colon [M]\mapsto 
  \( \theta_{\gamma_0}(M) : \dots : \theta_{\gamma_N}(M) \)
\end{equation}
  is a scheme-theoretic closed embedding.
\end{theorem}

\begin{proof}
In \eqref{eq:Thetagam}, one should interpret $\theta_{\gamma_i}(M)$
as the value at $[M]$ of the section $\theta_{\gamma_i}(a,b)$ of 
the ample line bundle $\lambda_U(a,b)$ on $\modspc{ss}{A}(a,b)$.
Thus we are simply describing a morphism to $\PP^N$ given by 
the linear system spanned by $N+1$ sections of a line bundle.

We use the identification in Proposition~\ref{prop:lamample} 
of sections of $\lambda_U(a,b)$ with semi-invariants
on $\Rep$, together with the fact from \cite{DW,SV} that 
such semi-invariants are spanned by determinantal ones,
to deduce that we can always choose a basis of sections of 
$\lambda_U(a,b)$ of the form $\theta_{\gamma_i}(a,b)$.

Hence the result follows by choosing
$U_0,U_1$ so that $\lambda_U(a,b)$ is very ample.
This is possible by the observation 
preceding Proposition~\ref{prop:lamample}.

Alternatively, the result can be proved
in a more generic, but less controlled way,
i.e. without the results of \secref{sub:theta-on-moduli}.
The construction of $\modspc{ss}{A}$ as a GIT quotient of the
representation space $\Rep$ means that it may be written as $\Proj(\cS)$,
where $\cS$ is the ring of semi-invariant functions on $\Rep$.
Thus, for some large $k$, there is a projective embedding determined 
by the linear system $\cS_k$, which has a basis of determinantal
semi-invariants, by \cite{DW,SV}, giving the required result.
\end{proof}

\begin{remark}
By the universal property described in Proposition~\ref{prop:descentA},
one may also interpret the morphism 
$\Theta_\gamma\colon \modspc{ss}{A}(a,b) \to \PP^N$ of \eqref{eq:Thetagam}
in a more functorial way as the unique morphism associated to
a natural transformation of functors 
$\Theta^\natural_\gamma\colon \modfun{ss}{A}(a,b)\to \funpt{\PP^N}$
by the fact that $\modspc{ss}{A}$ corepresents $\modfun{ss}{A}$
(cf. Definition~\ref{def:corep}).

The natural transformation $\Theta^\natural_\gamma$ 
is defined by essentially the same formula as \eqref{eq:Thetagam}, 
i.e. it takes $[M]$ in $\modfun{ss}{A}(S)$ to the element of $\funpt{\PP^N}(S)$
represented by the line bundle $\lambda_U(M)$ and the base-point free
linear system 
$\left\langle \theta_{\gamma_0}(M) , \dots , \theta_{\gamma_N}(M) \right\rangle$
(cf. \cite[Lecture~5]{Mu2} or \cite[Ch II, Th 7.1]{Ha}).

A similar remark applies to \eqref{eq:Thetadel} below.
\end{remark}

As a corollary of Theorem~\ref{thm:sep-Amod}, 
using essentially the adjunction in \eqref{eq:adj-on-hom},
we obtain a similar result for sheaves, 
with the usual more delicate conditions on the embedding.

\begin{theorem}
\label{thm:sep-OXmod}
  For any Hilbert polynomial $\hp$, we can find
  vector spaces $U_0,U_1$ satisfying \eqref{eq:UP}
  and finitely many maps
\[
  \delta_0,\dots,\delta_N \colon U_1\otimes \OX(-m) \to U_0\otimes \OX(-n)
\]
  such that the map 
\begin{equation}
\label{eq:Thetadel}
  \Theta_\delta \colon \modspc{ss}{X} \to \PP^N\colon [E] \mapsto 
  \( \theta_{\delta_0}(E) : \dots : \theta_{\delta_N}(E) \)
\end{equation}
 is a closed set-theoretic embedding.
 This embedding is scheme-theoretic in characteristic zero, while
 in characteristic $p$ it is scheme-theoretic on the stable locus.
\end{theorem}

\begin{proof}
We obtain the embedding, and its properties, 
by combining the embedding $\vphi\colon\modspc{ss}{X}\to \modspc{ss}{A}$
of Proposition~\ref{prop:embedding-final}
and the embedding $\Theta_\gamma\colon \modspc{ss}{A}\to \PP^N$
of Theorem~\ref{thm:sep-Amod}.
To see that $\Theta_\delta=\Theta_\gamma \circ \vphi$,
we need the observation from Proposition~\ref{prop:descentX}
that $\vphi^*\theta_{\gamma_i}=\theta_{\delta_i}$, for $\delta_i=\Phiadj(\gamma_i)$.
%
\end{proof}

\begin{remark}
\label{rem:theta-ring}
In characteristic zero, by considering the regularity of the ideal
sheaf of the embedding $\modspc{ss}{X}\subset \modspc{ss}{A}$,
we may deduce that, for sufficiently large $U_0,U_1$,
the restriction map 
\[
  \vphi^* \colon H^0(\modspc{ss}{A}(a,b),\lambda_U(a,b)) \to
  H^0(\modspc{ss}{X}(\hp),\lambda_U(\hp))
\]
is surjective.
Hence, as $H^0(\modspc{ss}{A},\lambda_U)$ is always spanned by 
determinantal theta functions $\theta_\gamma$,
we deduce that $H^0(\modspc{ss}{X},\lambda_U)$
is spanned by determinantal theta functions $\theta_\delta$,
for sufficiently large $U_0,U_1$.
\end{remark}

\subsection{Faltings' theta functions on curves} 
\label{sub:char-curve}

\newcommand{\g}{f}
\newcommand{\gp}{g}

A result of Faltings gives the following 
cohomological characterisation of
semistable bundles on a curve 
(see \cite{Fa} in characteristic zero 
and \cite{Se2} in arbitrary characteristic).

\begin{theorem}
\label{thm:cohochar-ssbundles-curves}
Let $X$ be a smooth projective curve. 
A bundle $E$ on $X$ is semistable if and only if
there exists a non-zero bundle $F$ such that
\begin{equation}
\label{eq:perpX}
 \Hom_X(F,E)=0=\Ext^1_X(F,E). 
\end{equation}
\end{theorem}

This condition \eqref{eq:perpX} may be interpreted 
as saying that $E$ and $F$ are
`perpendicular' in the derived category $\D(X)$,
in the sense that $\R\Hom_X(F,E)=0$.
Furthermore, this has the immediate numerical consequence that
\begin{equation}
  \label{eq:chi-def}
  \chi(F,E) \defeq \sum_{i\geq0} (-1)^i \dim \Ext_X^i(F,E) = 0. 
\end{equation}

Just as in \secref{sub:det-char-ss}, 
Theorem~\ref{thm:cohochar-ssbundles-curves} may also be interpreted as 
saying that certain determinantal theta functions 
detect the semistability of bundles on smooth curves.

To be precise about what this means, suppose that
$E$ is a family over $S$ of bundles on $X$
and $F$ is a bundle such that $\chi(F,E)=0$. 
Then $\R\Hom_X(F,E)$ is represented (locally over $S$) by a 
complex $d\colon \cK^0\to \cK^1$
of vector bundles, of the same rank,
such that, fibrewise at each $s\in S$,
$\ker d_s=\Hom_X(F,E_s)$ and $\coker d_s=\Ext^1(F,E_s)$.
There is then a well-defined line bundle $\lambda_F$ defined globally
on $S$, which is canonically isomorphic (locally) to 
$\det(\cK^0)^{-1}\otimes \det(\cK^1)$ and with a section $\theta_F$
canonically identified (locally) with $\det d$.
Note (see e.g. \cite[Lemma~2.5]{Se2}) that, 
if $r(F_1)=r(F_2)$ and $\det F_1=\det F_2$,
then $\lambda_{F_1}=\lambda_{F_2}$, so that ratios between
such theta functions $\theta_{F_1}$ and $\theta_{F_2}$
can meaningfully provide projective coordinates.

Indeed, Faltings \cite{Fa} shows that it is possible to find finitely many 
$F_0,\dots,F_N$ which detect all semistable bundles (of given rank $r$ and
degree $d$) and for which the morphism on the corresponding moduli space,
\[
  \Theta_F\colon \modspc{ss}{X}(r,d) \to \PP^N 
  \colon [E] \mapsto (\theta_{F_0}(E):\dots:\theta_{F_N}(E)),
\]
is the normalisation of its image, 
thereby giving an implicit construction of the moduli space.
Seshadri \cite[Remark~6.1]{Se2} asks how close this normalisation
is to being an isomorphism, or indeed, 
how close the theta functions $\theta_F$
come to spanning the space sections of the theta bundle $\lambda_F$ 
on $\modspc{ss}{X}$. 
Esteves \cite[Theorems~15,18]{Es} made progress by showing
that one can find a $\Theta_F$ which is a set-theoretic embedding
and which, in characteristic zero, is a scheme-theoretic embedding
on the stable locus $\modspc{s}{X}$.

Now, using Theorem~\ref{thm:sep-OXmod} and Remark~\ref{rem:theta-ring}, 
we can give a complete answer to Seshadri's question,
at least in characteristic zero.

First note that, by placing some reasonable restrictions 
on $E$ and $F$ it is possible to define theta functions globally.

\begin{lemma}
\label{lem:global-thetaF}
  Suppose that $E$ is a family over $S$ of $n$-regular sheaves
and that, for some $F$ with $\chi(F,E)=0$, there is a short exact sequence
\begin{equation}
  \label{eq:sesF}
   0\to F' \lra{\g} U\otimes \OX(-n) \lra{} F \to 0. 
\end{equation}
Then $\g^*\colon\Hom_X(U,E(n))\to\Hom_X(F',E)$ is
a map of vector bundles on $S$, of the same rank, and
$\theta_F=\det\g^*$.
\end{lemma}

\begin{proof}
 For any $s\in S$, apply the functor $\Hom_X(\?,E_s)$ to the short 
exact sequence \eqref{eq:sesF}. 
The resulting long exact sequence has just six terms, because $X$ is a
smooth curve. The fifth term $\Ext^1_X(U\otimes \OX(-n),E_s)$ vanishes
because $E_s$ is $n$-regular. 
Hence the sixth term $\Ext^1_X(F',E_s)$ also vanishes.

This vanishing means that $\Hom_X(F',E_s)$ is the fibre of a vector bundle
$\Hom_X(F',E)$ of rank $\chi(F',E)$
and that $\Hom_X(U,E_s(n))$ is the fibre of a vector bundle
$\Hom_X(U,E(n))$ of rank $\hp(n)\dim U$, which is equal to
$\chi(F',E)$, because $\chi(F,E)=0$.

Now, the remainder of the long exact sequence shows that the map
\[
  (\g^*)_s\colon \Hom_X(U,E_s(n))\lra{} \Hom_X(F',E_s)
\]
has kernel $\Hom_X(F,E_s)$ and cokernel $\Ext^1_X(F,E_s)$,
so that $\g^*$ represents $\R\Hom_X(F,E)$ (globally) 
and hence $\theta_F=\det\g^*$, as required.
\end{proof}

Using this we can show that the determinantal theta functions 
$\theta_\delta$ from \secref{sub:det-char-ss} are also theta
functions in the sense of Faltings.

\begin{proposition}
\label{prop:twothetas}
Let $\delta\colon U_1\otimes\OX(-m)\to U_0\otimes\OX(-n)$
be such that $\Hom_X(\delta,E_0)$ is invertible for some
bundle $E_0$ of rank $r$ and degree $d$
and let $F=\coker\delta$.

Then $\lambda_U\cong \lambda_F$ and 
$\theta_\delta=\theta_{F}$,
on any family of $n$-regular sheaves with rank $r$ and degree $d$.
In particular, for $E$ in such a family,
$\Hom_X(\delta,E)$ is invertible if and only if 
$\R\Hom_X(F,E)=0$. 
\end{proposition}

\begin{proof}
 We now have two short exact sequences
\begin{align}
    & 0\to F' \lra{\g} U_0\otimes \OX(-n) \lra{} F \to 0, \label{eq:ses0}\\
    & 0\to F'' \lra{} U_1\otimes \OX(-m) \lra{\gp} F' \to 0,\label{eq:ses1}
\end{align}
where $F'=\im\delta$, $F''=\ker\delta$ and $\delta=\g\gp$.

Note that, because $X$ is a smooth curve and hence its category of
coherent sheaves has global dimension 1,
we know that $\delta$, regarded as a complex $K_\bullet$,
is quasi-isomorphic to the direct sum of its cokernel and its (shifted)
kernel. 
Thus, $\Hom_X(\delta,E)$ is invertible 
if and only if $\R\Hom(K_\bullet,E)=0$ (see Remark~\ref{rem:dercat}),
which in turn is if and only if $\R\Hom(F,E)=0$ and $\R\Hom(F'',E)=0$.
Because we are assuming that this happens for one bundle $E_0$
of rank $r$ and degree $d$, this implies, in particular,
that $\chi(F,E)=0=\chi(F'',E)$, for any sheaf $E$
of the same rank and degree, by Riemann-Roch.
 
Now, we also observe that, for any $E$, 
we have the following factorisation of $\Hom_X(\delta,E)$, 
written here as $\delta^*$.
\begin{equation}\begin{gathered}
\label{diag:delta-curves}
  \xymatrix{
   \Hom_X(U_0,E(n)) \ar[d]^-{\g^*} \ar[rd]^-{\delta^*} & \\
   \Hom_X(F',E) \ar[r]^-{\gp^*}  & \Hom_X(U_1,E(m))
  }
\end{gathered}\end{equation}
The horizontal map $\gp^*$ is always injective.
Thus, to prove the equality of theta functions, 
what we need to show is that, when $E$ is $n$-regular, 
$\gp^*$ is an isomorphism, so that $\lambda_U\cong \lambda_F$ and 
\[
  \theta_\delta=\det \delta^*=\det \g^*=\theta_F,
\]
where the last equality is by Lemma~\ref{lem:global-thetaF}.

From the long exact sequence obtained by applying $\Hom_X(\?,E)$
to \eqref{eq:ses1} we see that it is sufficient to show that
$\Hom_X(F'',E)=0$.
We also see from the same long exact sequence that, 
when $E$ is $n$-regular, $\Ext^1_X(F'',E)=0$
and so the result follows, because we know that $\chi(F'',E)=0$.
\end{proof}

Note that the $F$ that occur here are necessarily vector bundles.

\begin{corollary}
\label{cor:Seshadri-answer}
  For given $r,d$, 
  there exists a finite set $F_0,\dots,F_N$ of vector bundles
  such that the map 
\[
  \Theta_F \colon \modspc{ss}{X}(r,d)\to \PP^N\colon [E] \mapsto 
  \( \theta_{F_0}(E) : \dots : \theta_{F_N}(E) \)
\]
 is a closed set-theoretic embedding.
 This embedding is scheme-theoretic in characteristic zero, while
 in characteristic $p$ it is scheme-theoretic on the stable locus.
\end{corollary}

\begin{proof}
  Immediate from Theorem~\ref{thm:sep-OXmod}
and Proposition~\ref{prop:twothetas}.
\end{proof}

Thus, in characteristic zero,
we see that Faltings' determinantal theta functions
can be used to give projective embeddings of the moduli spaces
of semistable bundles on a smooth curve.
Furthermore, by Remark~\ref{rem:theta-ring},
we have a positive answer to Seshadri's question:
the theta functions $\theta_F$ span the sections of line bundles $\lambda_U$
of sufficiently high degree.



\end{document}